\documentclass[10pt,twocolumn,amsmath,amssymb,pra,secnumarabic,%
    nofootinbib]{revtex4-1}

\usepackage{mathptmx,amsthm,enumitem,tikz,mathdots}
\usepackage[colorlinks=true,linkcolor=red!70!black,citecolor=green!70!black,
    urlcolor=magenta!70!black,backref]{hyperref}
\setlist[enumerate,1]{font=\bfseries,label=\arabic*.}
\newtheorem{theorem}{Theorem}[section]
\newtheorem{proposition}[theorem]{Proposition}
\newtheorem{lemma}[theorem]{Lemma}
\newtheorem{corollary}[theorem]{Corollary}
\newtheorem{conjecture}[theorem]{Conjecture}
\theoremstyle{remark}
\newtheorem*{remark}{Remark}

\makeatletter\def\@bibdataout@init{}\def\pre@bibdata{}\makeatother

\providecommand{\affiliation}{\address}
\providecommand{\acknowledgments}{\section*{Acknowledgments}}
\providecommand{\tensor}{\otimes}

\newcommand{\ER}{\mathsf{ER}}
\newcommand{\PR}{\mathsf{PR}}
\newcommand{\EXP}{\mathsf{EXP}}
\newcommand{\NP}{\mathsf{NP}}
\newcommand{\PSPACE}{\mathsf{PSPACE}}
\newcommand{\RE}{\mathsf{RE}}
\newcommand{\ccE}{\mathsf{E}}
\newcommand{\ccP}{\mathsf{P}}
\newcommand{\ccR}{\mathsf{R}}
\newcommand{\coRE}{\mathsf{coRE}}

\newcommand{\C}{\mathbb{C}}
\newcommand{\E}{\mathbb{E}}
\newcommand{\N}{\mathbb{N}}
\newcommand{\Q}{\mathbb{Q}}
\newcommand{\R}{\mathbb{R}}
\newcommand{\Z}{\mathbb{Z}}
\renewcommand{\H}{\mathbb{H}}

\newcommand{\GL}{\mathrm{GL}}
\newcommand{\ISO}{\mathrm{ISO}}
\newcommand{\Isom}{\mathrm{Isom}}
\newcommand{\PSL}{\mathrm{PSl}}
\newcommand{\SL}{\mathrm{SL}}

\newcommand{\SU}{\mathrm{SU}}
\newcommand{\Vol}{\mathrm{Vol}}
\newcommand{\bQ}{\bar{\Q}}
\newcommand{\del}{\partial}
\newcommand{\eps}{\epsilon}
\newcommand{\hN}{\hat{N}}
\newcommand{\hQ}{\hat{\Q}}
\newcommand{\mg}{\mathfrak{g}}
\newcommand{\no}{\mathrm{no}}
\newcommand{\poly}{\mathop{\mathrm{poly}}}
\newcommand{\tF}{\tilde{F}}
\newcommand{\tN}{\tilde{N}}
\newcommand{\te}{\tilde{e}}
\newcommand{\tv}{\tilde{v}}
\newcommand{\tw}{\tilde{w}}
\newcommand{\yes}{\mathrm{yes}}
\renewcommand{\mp}{\mathfrak{p}}
\renewcommand{\sharp}{\mathop{\#}}
\renewcommand{\tensor}{\otimes}

\newcommand{\congq}{\stackrel{?}{\cong}}

\newcommand{\Thm}[1]{Theorem~\ref{#1}}
\newcommand{\Lem}[1]{Lemma~\ref{#1}}
\newcommand{\Sec}[1]{Section~\ref{#1}}
\newcommand{\Prop}[1]{Proposition~\ref{#1}}
\newcommand{\Cor}[1]{Corollary~\ref{#1}}
\newcommand{\Conj}[1]{Conjecture~\ref{#1}}
\newcommand{\Fig}[1]{Figure~\ref{#1}}
\newcommand{\ie}{\emph{i.e.}}
\newcommand{\Ie}{\emph{I.e.}}

\newenvironment{eq}[1]{\begin{equation}\label{#1}}
    {\end{equation}\ignorespacesafterend}

\colorlet{lightblue}{white!85!blue}
\colorlet{darkred}{red!70!black}

\begin{document}

\title{Algorithmic homeomorphism of 3-manifolds as a corollary of
    geometrization}

\author{Greg Kuperberg}
\email{greg@math.ucdavis.edu}
\thanks{Partly supported by NSF grants CCF-1319245 and CCF-1716990.}
\affiliation{University of California, Davis}

\begin{abstract}
In this paper we prove two results, one semi-historical and the other new.
The semi-historical result, which goes back to Thurston and Riley, is
that the geometrization theorem implies that there is an algorithm for
the homeomorphism problem for closed, oriented, triangulated 3-manifolds.
We give a self-contained proof, with several variations at each stage,
that uses only the statement of the geometrization theorem, basic hyperbolic
geometry, and old results from combinatorial topology and computer science.
For this result, we do not rely on normal surface theory, methods from
geometric group theory, nor methods used to prove geometrization.

The new result is that the homeomorphism problem is elementary recursive,
\ie, that the computational complexity is bounded by a bounded tower of
exponentials.  This result relies on normal surface theory, Mostow rigidity,
and bounds on the computational complexity of solving algebraic equations.
\end{abstract}

\maketitle

\section{Introduction}

In this paper, we will prove the following two theorems.

\begin{theorem}[After Thurston \cite{Thurston:kleinian}] Suppose that
$M_1$ and $M_2$ are two finite, simplicial complexes that represent closed,
oriented 3-manifolds.  Then, as a corollary of the geometrization theorem, it
is recursive to determine if there is an orientation-preserving homeomorphism
$M_1 \cong M_2$.
\label{th:main} \end{theorem}

\begin{theorem} The oriented homeomorphism problem for closed, oriented
3-manifolds is elementary recursive.
\label{th:ermain} \end{theorem}

\Thm{th:main} implies that the geometrization theorem is a classification of
closed, oriented 3-manifolds by the standard of computer science, where the
term \emph{recursive} used here means the same thing as \emph{decidable}
or \emph{computable}.  Geometrization intuitively presents itself as a
classification of closed 3-manifolds, or at least a big step towards one.
However, the question of what counts as a ``classification" in mathematics is
generally not rigorous, even though it is typically a debate over rigorous
results.  The computability interpretation is thus important because it is
rigorous, even though it is not by any means the only important standard of
classification.  (For instance, the set of twin primes is recursive, but they
remain unclassified in the sense that it is not even proven that there are
infinitely many.)  Note that Thurston himself \cite[Sec.~3]{Thurston:kleinian}
seriously addressed the relation between geometrization and computability.

We argue that \Thm{th:main} should largely be credited to Riley and
Thurston from the 1970s, even though they did not publish a complete proof.
(See \Sec{ss:history} for more details.)  To support this interpretation,
we will prove \Thm{th:main} directly using hyperbolic geometry; and using
other background results concerning computability and triangulations of
manifolds that seem standard and germane.  The most important results of
the latter type are the Tarski-Seidenberg theorem (\Thm{th:tarski}) that real
algebraic equations can be solved recursively, Kantorovich's theorem on
convergence of Newton's method (\Thm{th:nkn}), and the stellar and bistellar
move theorems of Alexander, Newman, and Pachner (Theorems~\ref{th:stellar}
and \ref{th:bistellar}).  Despite this restriction on methods, we give
more than one argument for each of several stages of the proof.

In the intervening years, Jaco-Tollefson, Manning, Scott-Short,
and others have published proofs of major parts of \Thm{th:main}
\cite{JT:complete,Manning:algo,SS:problem,AFW:decision}.  These approaches
have various new ideas and implications, which is in keeping with
Thurston's philosophy concerning the nature of progress in mathematics
\cite{Thurston:progress}.   Even so, the status of \Thm{th:main} has remained
unsettled.  At one extreme, it has been interpreted as a folklore theorem
and therefore standard knowledge, even if the proof is not elementary.
At the other extreme, it has been interpreted as still an open problem.
In the middle, one could argue that the published partial results piece
together to make an entire proof.  The problem with the middle position
is that the total structure of an arbitrary closed, oriented 3-manifold
is somewhat complicated; see Theorems~\ref{th:factor}, \ref{th:jsj}, and
\ref{th:geom}.  So, one purpose of our proof of \Thm{th:main} is to give
a complete proof in one paper, as requested by Aschenbrenner-Friedl-Wilton
\cite{AFW:decision}.

The intervening results also typically use either normal surface
theory \cite{Kneser:flachen,Haken:normal} or geometric group theory
\cite{Gromov:hyperbolic,Sela:isomorphism,ECHLSPT:word}.  While these methods
certainly work, they arguably overshoot \Thm{th:main}.  Both theories are
highly non-trivial in their own right, and they continued to be developed
after the geometrization conjecture was stated.  In particular the key
results of Jaco-Tollefson \cite{JT:complete} and Sela \cite{Sela:isomorphism}
came later.  Sela's theorem applies to Gromov-hyperbolic groups, which are
vastly more general than the Kleinian groups that arise in \Thm{th:main}.
Meanwhile Haken and Jaco-Tollefson prove sharper results than strictly
necessary for their components of \Thm{th:main}; namely, they establish
algorithms with quantitative bounds on execution time.  This brings us
to \Thm{th:ermain}.

In \Thm{th:ermain}, an algorithm is \emph{elementary recursive} if its
execution time is bounded by a bounded tower of exponentials; for instance,
time $O(2^{2^n})$.  (See \Sec{ss:errec}.)  In contrast with \Thm{th:main},
the proof of \Thm{th:ermain} does use normal surface theory, as well as Mostow
rigidity, and improved bounds on the computational complexity of solving
algebraic equations \cite{Grigoriev:algebra}.  The connected-sum and JSJ
decomposition stages of \Thm{th:ermain} were partly known.  For instance,
using similar methods, Mijatovi\'c \cite{Mijatovic:seifert,Mijatovic:s3}
established an elementary recursive bounds on the number of Pachner moves
needed to standardize either $S^3$ or a Seifert-fibered space with boundary.

The hyperbolic case of \Thm{th:ermain} is new.  By contrast, Mijatovi\'c
also established a primitive recursive bound on the number of Pachner moves
needed to equate two hyperbolic, fiber-free, Haken 3-manifolds.   However,
primitive recursive is significantly weaker than elementary recursive;
the Haken condition is also a significant restriction.  \Thm{th:ermain}
also has the advantage of combining a mixed set of methods to handle the
full generality of closed, oriented 3-manifolds.

\begin{remark} We leave the non-orientable versions of Theorems \ref{th:main}
and \ref{th:ermain} for future work.  This case includes new details such as
3-manifolds with essential, two-sided projective planes and Klein bottles.
A more thorough result would also handle compact 3-manifolds with boundary.
\end{remark}

\subsection{History and discussion}
\label{ss:history}

As already mentioned, the geometrization conjecture has often been
interpreted as a classification of closed 3-manifolds, and 
computability is one candidate standard of what it means to classify
mathematical objects.  In Thurston's famous survey of his results
in the AMS Bulletin \cite[Sec.~3]{Thurston:kleinian}, he says:
\begin{quote} Riley's work makes it clear that there is a rigorous, but
not generally practical, algorithm for computing hyperbolic structures.
\end{quote}
Thurston then sketches an algorithm which is similar Manning's construction
\cite{Manning:algo} in some ways and to our arguments in other ways.
This passage, and some other aspects of the Bulletin article, support the
conclusion that Thurston anticipated not only the statement of \Thm{th:main},
but also its proof.  The author also discussed \Thm{th:main} in personal
communication with Thurston in the late 1990s.

At first glance, an algorithm that can only find a hyperbolic structure on
a 3-manifold is both less general and weaker than \Thm{th:main}. It is less
general because a 3-manifold may also have non-hyperbolic components; it is
weaker because the homeomorphism problem for two hyperbolic manifolds $M_1$
and $M_2$ takes more work than just finding their hyperbolic structures.
However, in the theorem that it is recursive to geometrize a 3-manifold
(\Thm{th:geomrec}), the hyperbolic pieces (\Lem{l:hyprec}) are the hardest
part. The geometric structure of the other pieces and the data to glue the
pieces together are complicated to describe carefully, but the proof that
this data is computable only requires moves on triangulations (\Cor{c:plre}),
the principle that nested infinite loops can be combined into a single
infinite loop (\Prop{p:gsearch}), and some facts about Seifert fibrations
(\Lem{l:seifrec} and \Thm{th:waldhausenb}).

In the second stage, the homeomorphism problem for hyperbolic 3-manifolds
(\Thm{th:hyphom}) reduces to calculating isometries by Mostow rigidity, and
a typical algorithm for this is similar to one for computing a hyperbolic
structure.  The homeomorphism problem for Seifert-fibered components and
glued combinations of components (\Sec{s:homrec}) requires little more than
the ideas of Waldhausen \cite{Waldhausen:gruppen} in his classification
of graph manifolds.

Later in the Bulletin article \cite[Sec.~6]{Thurston:kleinian}, Thurston
gives a list of open problems and projects, including: 
\begin{quote}21. Develop a computer program to calculate hyperbolic
structures on 3-manifolds.
\end{quote}
Jeff Weeks' SnapPea \cite{W:snappea} (now SnapPy \cite{W:snappy}), which
was originally written in the 1980s, met this challenge.  It is fast
and reliable in practice, it can also compute the isometries between
two hyperbolic 3-manifolds, and it has been supremely useful for a lot of
research in 3-manifold topology.  SnapPea also supports the belief that the
homeomorphism problem follows from geometrization, given its spectacular
record in practice.  However, its specific algorithms are not rigorous.
SnapPea uses ideal triangulations of cusped 3-manifolds, together with
Dehn fillings to make spun triangulations of closed 3-manifolds; it is
only conjectured that such a structure always exists.  SnapPea also uses
non-rigorous methods to find suitable triangulations.   In particular
it uses limited-precision floating point arithmetic; it has no rigorous
model of necessary precision as a function of geometric complexity.
(Note that SnapPy has been extended to rigorously certify an answer,
when SnapPea finds one.)   In contrast to the SnapPea data structure, and
other reasons that ideal and spun triangulations are important, we will use
triangulations with semi-ideal and finite tetrahedra to prove \Thm{th:main}
(see \Sec{sss:hypcomp}).

To start the rigorous discussion of computable classification and the
homeomorphism problem, we can say that closed 3-manifolds are classified
if we can:
\begin{enumerate}
\item specify every closed 3-manifold by a finite data structure;
\item algorithmically generate a standard list of closed 3-manifolds without
repetition; and
\item given any 3-manifold $M$, algorithmically identify the standard
manifold $M'$ such that $M \cong M'$.
\end{enumerate}
For closed 3-manifolds, condition 1 is addressed by the fact that every
3-manifold has a unique smooth structure and a unique PL structure.  As a
result, we can describe a closed 3-manifold as a finite simplicial complex.
Unlike in higher dimensions, it is easy to check whether a simplicial complex
is a 3-manifold (\Sec{s:triang}).  Conditions 2 and 3 are equivalent to
an algorithm to determine whether two closed 3-manifolds $M_1$ and $M_2$
are homeomorphic by the following simple argument.  (Haken calls this
argument the ``cheapological trick" \cite[Sec.~4]{Waldhausen:recent}.)
In one direction, if both conditions 2 and 3 are satisfied, then Condition
3 immediately implies a homeomorphism algorithm.  In the other direction,
given a homeomorphism algorithm, we can make can lexicographically order all
descriptions of all closed 3-manifolds according to condition 1, and then
list only those examples that are not homeomorphic to any earlier example.
This satisfies condition 2.  Then given a description of a closed 3-manifold
$M$, we can search the list in order to find the standard $M' \cong M$
to satisfy condition 3.  (Haken calls this argument the ``cheapological
trick" \cite[Sec.~4]{Waldhausen:recent}.  Arguably it is not cheap after
all, since it is substantially similar to actual constructions of
tables of knots and 3-manifolds.)

As mentioned, Manning \cite{Manning:algo} and Scott and Short
\cite{SS:problem} give partial results toward \Thm{th:main}, but
they use more recent tools.  In particular, Manning uses Sela's algorithm
\cite{Sela:isomorphism} for the isomorphism problem for word-hyperbolic
groups, while Scott and Short use the theory of biautomatic groups
\cite{ECHLSPT:word}.

Both Short-Scott and Aschenbrenner-Friedl-Wilton
\cite[Sec.~2.1]{AFW:decision} mention a particular subtlety in approaches to
\Thm{th:main} that are based on analyzing the fundamental group $\pi_1(M)$
or the fundamental groups of its components.  Namely, $\pi_1(M)$ is
insensitive to the orientation of $M$.  Worse, if
\[ M \cong W_1 \sharp W_2 \sharp \dots \sharp W_n \]
is a decomposition of $M$ into prime summands, then the orientation of each
summand $W_k$ can be chosen separately without changing $\pi_1(M)$.  Or the
summands can be lens spaces; two lens spaces can have the same fundamental
group without even being unoriented homeomorphic.  We surmount this subtlety
by modelling all 3-manifolds and their components with triangulations that
are decorated with orientations; see \Sec{ss:compgeom}.

\acknowledgments

The author would like to thank David Futer, Bus Jaco, Misha Kapovich,
Jason Manning, Peter Scott, and Henry Wilton for very useful discussions.
The author would also like to thank the anonymous referee for a careful
reading that led to significant revisions.

\section{Computability}

We recommend Arora-Barak \cite{AB:modern} and the Complexity Zoo \cite{W:zoo}
for modern introductions to models of computation and complexity classes.

\subsection{Recursive and recursively enumerable problems}

Let $A$ be a finite alphabet and let $A^*$ be the set of all
finite words over that alphabet.  A \emph{decision problem} is a function
\[ d:A^* \to \{\yes,\no\}. \]
A \emph{function problem} is likewise a function $f:A^* \to A^*$,
which can be multivalued.  The input space $A^*$ is equivalent to
many other types of input by some suitable encoding:  Finite sequences of
strings, finite simplicial complexes, etc.

A decision problem or a function problem can be a \emph{promise
problem}, meaning that it is defined only on some subset of inputs $P
\subseteq A^*$ which is called a \emph{promise}.  Whether two closed
$n$-manifolds are PL homeomorphic is an example of a promise decision
problem:  The input consists of two simplicial complexes that are promised
to be manifolds; then the yes/no decision is whether they are homeomorphic.
(But see \Prop{p:manifold}.)

A \emph{decision algorithm} is a mathematical computer program, which can
be modelled by a Turing machine (or some equivalent model of computation),
that takes some input $x \in A^*$ and can do one of three things:
(1) Terminate with the answer ``yes", (2) terminate with the answer ``no",
or (3) continue in an infinite loop.  Similarly, a function algorithm can
terminate and report an output $y \in A^*$, or it can continue in an
infinite loop.  Given a multivalued function $f$, then a function algorithm
is only required to calculate one of the values of $f(x)$ on input $x$.

A \emph{complexity class} or \emph{computability class} is some set
of decision or function problems, typically defined by the existence of
algorithms of some kind.  For example, a decision problem $d$ or a function
problem $f$ is \emph{recursive} (or \emph{computable} or \emph{decidable})
if it is computed by an algorithm that always terminates.  By definition,
the complexity class $\ccR$ is the set of recursive decision problems.
By abuse of notation, $\ccR$ can also denote the set of recursive,
promise decision problems; or the set of recursive function problems,
with or without a promise.  The following proposition is elementary.

\begin{proposition} If $d$ is a recursive promise problem,
and if the promise itself is recursive, then $d$ is a recursive
non-promise problem if we let $d(x) = \no$ when the promise is not 
satisfied.
\label{p:promise} \end{proposition}

The complexity class $\RE$ is the set of \emph{recursively enumerable}
decision problems.   These are problems with an algorithm that terminates
with ``yes" when the answer is yes; but if the answer is ``no", the
algorithm might not terminate.  The complexity class $\coRE$ is defined
in the same way as $\RE$, but with yes and no switched.  
We review the following standard propositions and theorems.

\begin{proposition} A non-promise decision problem $d$ is in $\RE$ if and
only if there is an algorithm that lists all solutions to $d(x) = \yes$
without repetition.
\label{p:list} \end{proposition}

\Prop{p:list} justifies the name ``recursively enumerable" for the class
$\RE$.  (Note that the solutions can be listed \emph{in order} if and only
if $d \in \ccR$.)  The proof is left as an exercise.  Also, in the spirit
of \Prop{p:list}, a decision problem $d$ can be identified with the set
of solutions to $d(x) = \yes$; in this way we can call a set recursive,
recursively enumerable, etc.

\begin{proposition} $\ccR = \RE \cap \coRE$.
\label{p:recore} \end{proposition}

The proof of \Prop{p:recore} is elementary but important: Given separate
$\RE$ algorithms for both the ``yes" and ``no" answers, we can simply run
them in parallel; one of them will finish.  The proposition and its proof
reveal the important point that a recursive algorithm might come with no
bound whatsoever on its execution time.

\begin{theorem}[Turing] The halting problem is in $\RE$ but not in $\ccR$.
In particular, $\RE \ne \ccR$.  
\end{theorem}

Informally, the \emph{halting problem} is the question of whether a given
algorithm with a given input terminates.   Let $h$ be the halting decision
problem, where the input $x$ in the value $h(x)$ is an encoding of an
algorithm and its input (or, traditionally, an encoding of a Turing machine).
It is easy to show that $h$ is \emph{$\RE$-complete} in the following sense:
Given a problem $d \in \RE$, there is a recursive function $f$ such that
$d(x) = h(f(x))$ for any input $x$.  Any other problem in $\RE$ with this
same property can also be called $\RE$-complete, or \emph{halting complete}.

\begin{proposition} Let $G$ be a graph structure on $A^*$.  If the edge
set of $G$ is recursively enumerable, then so is the set of pairs $(x,y)$,
where $x$ and $y$ are vertices in the same connected component of $G$.
\label{p:gsearch} \end{proposition}

\Prop{p:gsearch} is important for recursively enumerable infinite searches.
The interpretation of the proposition, which is conveyed by the proof,
is that nested infinite loops can be reorganized into a single infinite loop.

\begin{proof} By \Prop{p:list}, we can model a recursively enumerable set
by an algorithm that lists its elements.  The proposition
states that the elements can be listed without repetition; but this
is optional, since we can store all of the elements already listed
and omit duplicates.

We use a recursive bijection $f$ between the natural numbers $\N$ and
$\N^*$, the set of finite sequences of elements of $\N$.   We can express
any element of $\N^*$ uniquely in a finite alphabet that consists of the
ten digits and the comma symbol.   We can then list of all of these strings
first by length, and then in lexicographic order for each fixed length.
We can then let $f(n)$ be the $n$th listed string.

We can now convert the value $f(n)$ to a finite path $(x_0,\ldots,x_k)$
in the graph $G$, in such a way that every finite path is realized.  If
\[ f(n) = (n_0,\ldots,n_k), \]
then we let $x_0$ be the $n_0$th string in $A^*$. For each $j > 0$,
we let $x_j$ be the $n_j$th neighbor of $x_{j-1}$.  In order to find the
$n_j$th neighbor of $x_{j-1}$, we list the elements of the edge set of $G$
until the edge $(x_{j-1},x_j)$ arises as the $j$th edge from $x_{j-1}$.
There is the technicality that $x_{j-1}$ might not have an $n_j$th
neighbor if it only has finitely many neighbors.  To avoid this problem,
we intersperse trivial edges of the form $(x,x)$ infinitely many times,
for every string $x \in A^*$, along with the non-trivial edges of $G$.

Since the algorithm finds every finite path in $G$, it finds every pair
of vertices $x$ and $y$ in the same connected component.  Thus, the set
of such pairs is recursively enumerable.
\end{proof}

\subsection{Elementary recursive problems}
\label{ss:errec}

As mentioned after \Prop{p:recore}, a recursive algorithm need not have
any explicit upper bound on its execution time, beyond the tautological
bound that running it is a way to calculate how long it runs.  This
motivates smaller complexity classes that are defined by explicit bounds.
The most common notation for a bound on the execution time of an algorithm
is asymptotic notation as a function of the input length $n = |x|$ to a
decision problem $d(x)$.  For example, we could ask for a polynomial-time
algorithm, by definition one that runs in time $O(n^k)$ for some fixed $k$.

We have two reasons to consider a fairly generous bound in this paper.
First, the recursive class $\ccR$ is unfathomably generous, so any explicit
bound can be considered a major improvement.   Second, the computational
complexity of a problem or algorithm depends somewhat on the specific
computational model, but certain relatively generous complexity bounds
are substantially model-independent.

We consider a traditional Turing machine first.
By (informal) definition, a \emph{Turing machine} is a finite-state ``head"
with an infinite linear memory tape, and deterministic dynamical behavior.
We say that an algorithm is \emph{elementary recursive} if it runs in time
\[ O\left(\underbrace{2^{2^{\iddots^{2^n}}}}_k\right) \]
for some constant $k$.  We call the corresponding complexity class $\ER$.
By abuse of terminology, we use $\ER$ to refer to both decision problems
and function problems, and to numerical bounds.  (Note that if $f(n)$ can
be computed in $\ER$, then a running time bound of $O(f(n))$ is itself a
subclass of $\ER$.)

Without reviewing rigorous definitions, we list some of the many
variations in the computational model that do not affect the class $\ER$
in the following proposition.   The proposition is not really needed in
this paper except to motivate $\ER$ as an important complexity class.
The only tacit dependence is that a random access machine is somewhat
closer to both intuitive descriptions of algorithms and actual computers
than a Turing machine with a linear tape is.

\begin{proposition} Each of the following four computational
models is the same as standard $\ER$.
\begin{enumerate}
\item A Turing machine with an $n$-dimensional tape, or a random access
tape addressable by a separate address tape, with an elementary recursive
bound on computation time.
\item A Turing machine restricted to an elementary
recursive bound on computational space and unrestricted computation time.
\item A randomized Turing machine whose answers are probably
correct, with an elementary recursive bound on computation time.
\item A quantum Turing machine that can compute in quantum
superposition, with an elementary recursive bound on computation time.
\end{enumerate}
\label{p:er} \end{proposition}

\begin{proof} Instead of a self-contained proof, we justify each case of
the proposition with specific references to Arora-Barak \cite{AB:modern}.
\begin{enumerate}
\item This follows from Exercises 1.7 and 1.9 in Arora-Barak.
\item This follows from Theorem 4.2 in Arora-Barak.
\item This reduces to case 4 by the proof of Corollary 10.11 in Arora-Barak.
\item This reduces to case 2 by the proof of Theorem 10.23 in Arora-Barak.
\qedhere
\end{enumerate}
\end{proof}

\begin{remark} An elementary recursive bound is also a major improvement
over another bound that is popular in logic and computer science: primitive
recursive.  An algorithm is \emph{primitive recursive} if it runs in time
$O(n[k]b)$ for some fixed $k$ and $b$, where the $k$th operation $a[k]b$
is defined inductively as follows:
\begin{gather*}
a[1]b = a+b \qquad a[2]b = ab \qquad a[3]b = a^b \\
a[k+1]b = \underbrace{a[k](a[k](\cdots (a[k]a)\cdots))}_b.
\end{gather*}
For example, the operation $a[4]b$, which is called \emph{tetration}, is
defined as a tower of exponentials of height $b$.  The primitive recursive
complexity class is denoted $\PR$.  We can organize $\PR$ into a complexity
hierarchy by defining $\ccE_k$ to be the set of functions computable in
time $O(n[k+1]b)$ for some fixed $b$.  Then $\ccE_2 = \ccP$, $\ccE_3 =
\ER$, $\ccE_4$ consists of complexity bounds which are bounded towers of
tetrations, etc.
\end{remark}

\subsection{Computable numbers}

A \emph{computable real number} $r \in \R$ is a real number with a
computable sequence of bounding rational intervals.   In other words,
there is an algorithm that generates rational numbers $a_n,b_n \in \Q$
such that $x \in [a_n,b_n]$ and $b_n-a_n \to 0$.  Many standard algorithms
from numerical analysis, including field operations, integration of
continuous functions, Newton's method, etc., have the property that if
the input consists of computable numbers, then so does the output. One
main limitation of computable real numbers is that inequality tests such
as $a > b$ or $a \ne b$ are only recursively enumerable, not recursive.
In other words, if $a \ne b$, then there is an algorithm to eventually
confirm this fact and say which one is greater; but there is no
terminating algorithm that always confirms that $a=b$.

We can avoid this shortcoming of the field of computable numbers by passing
to a smaller subfield where equality is also recursive.  In particular,
we will use $\hQ = \R \cap \bQ$, the real algebraic closure of the rational
numbers $\Q$, which has this property.

\begin{theorem} There is an encoding of the elements of $\hQ$ such
that field operations, order relations, and conversion to
computable real numbers are all recursive.
\label{th:realalg} \end{theorem}

One encoding of a real algebraic number $x$ that can be used to prove
\Thm{th:realalg} is to describe it by a minimal polynomial together with an
isolating interval $x \in [a,b]$ with rational endpoints to distinguish $x$
from its Galois conjugates.   Note that the isolating interval may be made
arbitrarily small since algebraic numbers are computable, for instance by
Newton's method.  Note also that a computable encoding of elements of $\hQ$
yields a computable encoding of elements of $\bQ \subseteq \C$ the field
of all algebraic numbers.

\begin{remark} The field of real algebraic numbers together with reliable
equality testing is implemented in \texttt{Sage} \cite{W:sage}.
\end{remark} 

\begin{theorem}[Tarski-Seidenberg \cite{Tarski:method,Seidenberg:method}]
It is recursive to determine whether there is a solution to a finite
list of polynomial equalities and inequalities with coefficients in $\hQ$
in finitely many variables; or to find a solution.
\label{th:tarski} \end{theorem}

Actually, Tarski and Seidenberg proved the stronger result that it is
recursive to decide any assertion over $\R$ expressed with polynomial
relations and first-order quantifiers.

\section{Triangulations of manifolds and moves}
\label{s:triang}

In this section, we will analyze the form of the input to \Thm{th:main}.
We will show that given simplicial complexes as input to the homeomorphism
problem, we can first confirm that they are 3-manifolds.   (It is also easy
to confirm that input strings actually represent simplicial complexes,
in some convenient data type.)   Thus \Prop{p:promise} applies: we can
view the homeomorphism problem as a non-promise problem.  Actually,
\Prop{p:manifold} below is overkill for this purpose, since it is much
harder in dimension $n=4$ than in dimension $n=3$.

We then discuss moves between triangulations of a manifold, mainly to
establish \Cor{c:plre}.  In light of \Prop{p:recore}, \Cor{c:plre} is an
easy half of \Thm{th:main}, one that holds in any dimension $n$.

\begin{proposition} If $\Theta$ is a finite simplicial complex of dimension
$n \le 4$, then it is recursive to determine whether it is a closed PL
$n$-manifold, and whether or not it is orientable.
\label{p:manifold} \end{proposition}

\begin{proof} The proof is partly by induction on dimension $n$.
The result is trivial if $n=0$, where we need only check that $\Theta$
is a single point.  Otherwise, we first check that $\Theta$ is connected,
and we must check that the link $\Lambda$ of every vertex is both a
closed $(n-1)$-manifold and a PL $n$-sphere.  The former condition is
the inductive step.   The latter condition requires an algorithm to
recognize an $(n-1)$-sphere.  If $\Lambda$ is a closed 1-manifold, then
it is immediately a 1-sphere, \ie, a circle.  If $\Lambda$ is a closed
2-manifold, then we can compute its Euler characteristic.  If $\Lambda$ is
a closed 3-manifold, then \Thm{th:main} implies that it is recursive to
determine if $\Lambda$ is a 3-sphere, although this result was obtained
without geometrization by Rubinstein and Thompson (\Thm{th:ersphere})
\cite{Rubinstein:recog,Thompson:recog}.

We can check that $\Theta$ is orientable (and orient it) algorithmically by
computing its simplicial homology.
\end{proof}

The stellar and bistellar subdivision theorems establish that every two
triangulations of a compact $n$-manifold, in particular a compact 3-manifold,
are connected by a finite sequence of explicit moves.  See Lickorish
\cite{Lickorish:complexes} for a modern treatment and a historical review.

\begin{theorem}[Alexander-Newman] If two finite simplicial complexes
$\Theta_1$ and $\Theta_2$ are PL equivalent, then they are connected by
a sequence of stellar subdivision moves and their inverses.
\label{th:stellar} \end{theorem}

Briefly, a \emph{stellar} move in a simplicial complex $\Theta$ consists of
replacing the star $\mathrm{st}(\Delta)$ of some simplex $\Delta$ in $\Theta$
with a cone over the subcomplex of simplices in $\mathrm{st}(\Delta)$
that do not contain $\Delta$.   Equivalently, the apex $v$ of this cone is
placed in the interior of $\Delta$, and all simplices that contain $\Delta$
are subdivided to support the new vertex $v$.

\begin{theorem}[Pachner] If $\Theta_1$ and $\Theta_2$ are two triangulations
of a compact, PL manifold $M$, then they are connected by bistellar moves.
\label{th:bistellar} \end{theorem}

A \emph{bistellar} move of a triangulation of an $n$-manifold $M$
consists of a stellation followed an inverse stellation at the same vertex.
Equivalently, two triangulations of $M$ differ by a bistellar move when there
is a minimal cobordism between them consisting of a single $(n+1)$-simplex.
In particular, a shellable triangulation of $M \times I$ yields a sequence
of bistellar moves.

Lickorish points out that Newman conjectured and partially proved
\Thm{th:bistellar} in an earlier paper, before he and Alexander separately
proved \Thm{th:stellar}.  Bistellar moves are also called Pachner moves,
although arguably they should be called Newman-Pachner moves.

\Thm{th:bistellar} also holds for ideal or semi-ideal triangulations of
a compact 3-manifold with torus boundary components.  (In other words,
it holds for a 3-dimensional pseudomanifold with singular points with
torus links, which are the ideal vertices.)

Theorems \ref{th:stellar} and \ref{th:bistellar} each have the following
corollary.

\begin{corollary} The PL homeomorphism problem for compact PL $n$-manifolds
is in $\RE$.
\label{c:plre} \end{corollary}

\begin{remark} There is a proof of \Cor{c:plre} that works directly from
the definition of PL equivalence without using Theorem~\ref{th:stellar}
or \ref{th:bistellar}, nor even \Prop{p:gsearch}.  For each $n$, choose
a linear embedding of an $n$-simplex $\Delta^n \subseteq \R^n$.  Then in
general a \emph{geometric refinement} is a simplicial complex $\Theta$
with a homeomorphism onto $\Delta$ which is affine-linear on each simplex
of $\Theta$.  Likewise a refinement of a simplicial complex $\Theta_1$
is another simplicial complex $\Theta_2$ with a homeomorphism $f:\Theta_2
\to \Theta_1$, such that $f$ yields a geometric refinement of each simplex
of $\Theta_1$.  By definition, two complexes $\Theta_1$ and $\Theta_2$
are PL equivalent if they share a refinement $\Theta_3$.  Now, we can
let each $\Delta^n$ have rational vertices (\ie, vertices in $\Q^n$), and
after that we can perturb any geometric refinement so that its vertices
are all rational.  The set of rational mutual refinements of two finite
complexes $\Theta_1$ and $\Theta_2$ is recursive by direct verification.
(In other words, given a simplicial complex $\Theta_3$, and given rational
target positions for its vertices in both $\Theta_1$ and $\Theta_2$, we
can algorithmically check whether this data yields a mutual refinement.)
Therefore the question of whether there exists a mutual refinement is
directly recursively enumerable.
\end{remark}

\begin{proposition} If $\Theta_1$ is a finite simplicial complex with
$n_1$ simplices (of arbitrary dimension) and $n_2 \ge n_1$, then it is
recursive to produce a complete list of geometric subdivisions $\Theta_2$
of $\Theta_1$ with $n_2$ simplices.
\label{p:refine} \end{proposition}

\begin{proof} There are only finitely many simplicial complexes $\Theta_2$
with $n_2$ simplices, and they can be generated recursively.  For each
candidate for $\Theta_2$, there are only finitely many combinatorial
choices for a function from the simplices of $\Theta_2$ to the simplices
of $\Theta_1$.  For each such choice, we can first check that the simplices
of $\Theta_2$ that land in a $k$-simplex $\Delta \in \Theta_1$ support a
simplicial cycle that represents the fundamental class in $H^k(\Delta,\del
\Delta)$.  We solve for each such cycle for all $\Delta$ (where each must
be unique if $\Theta_2$ indeed subdivides $\Theta_1$).  Then the constraint
that each simplex of $\Theta_2$ must be positively oriented in $\Theta_1$
yields we obtain algebraic inequalities for the positions of all vertices.
We can then apply \Thm{th:tarski} to see if there is a solution for those
positions.
\end{proof}

We conclude this section with the following theorem which combines
results of P.S. Novikov, Boone, Adian, Rabin, Markov, and S.P. Novikov
\cite{Poonen:sampler}.

\begin{theorem}[NBARMN] The isomorphism problem for finitely presented groups,
the PL homeomorphism problem for 4-manifolds, and the recognition of $S^n$
among PL $n$-manifolds for each $n \ge 5$ are all halting-complete.
\label{th:complete} \end{theorem}

It is not known whether either topological or PL recognition of $S^4$
is recursive.

\begin{remark} The homeomorphism problem for PL $n$-manifolds in
\Thm{th:complete}, or even recognition of $S^n$, needs to be handled with
some care, for several reasons.  First, because recognizing whether the input
is a PL $n$-manifold is (by \Thm{th:complete}!) an uncomputable promise when
$n \ge 6$.  Second, because there are closed manifolds that are homeomorphic
but not PL homeomorphic \cite{KS:hauptvermutung}.  Third, because there
are simplicial complexes that are not PL $n$-manifolds at all, but that
are homeomorphic to $S^n$, for each $n \ge 5$ \cite{Edwards:suspensions}.
The proof of \Thm{th:complete} dispenses with all of these concerns as
follows.  Given an input $x$ to the halting problem $h(x)$ and an integer
$n \ge 4$, there is an algorithm that constructs an $n$-manifold $M(x)$
such that:
\begin{enumerate}
\item $M(x)$ is manifestly a closed PL manifold.
\item $M(x)$ is PL homeomorphic to $S^n$ when $n \ge 5$, or to a
connected sum of copies of $S^2 \times S^2$ when $n=4$, if and only if
$M(x)$ is simply connected.
\item $M(x)$ is simply connected if and only if $h(x) = \yes$.
\end{enumerate}
\end{remark}

\begin{remark} By contrast with \Thm{th:complete}, the PL homeomorphism
problem for simply connected $n$-manifolds with $n \ge 5$ is recursive
\cite{NW:aspects}.
\end{remark}

\section{Some notation}
\label{s:notation}

We summarize some notation for specific topological spaces, beyond the
most standard notation that $S^n$ is an $n$-sphere, $D^n$ is an $n$-disk,
$P^n$ is real projective $n$-space, and $I = D^1$ is an interval.

We let $X \ltimes Y$ denote a fiber bundle with base $X$ and fiber $Y$.
Although the notation $X \tilde{\times} Y$ is reasonably standard for a
twisted bundle, we prefer to write $X \ltimes Y$, for two reasons.  First,
because the notation specifies which side is the fiber; we can write $X
\ltimes Y \cong Y \rtimes X$.  Second, because a fiber bundle is analogous
to a semidirect product in group theory.

We review Seifert's description of oriented Seifert-fibered spaces
\cite{Seifert:fibered}.  If $F$ is a compact surface which may
or may not be orientable, then there is a unique, canonically oriented
$I$-bundle $F \ltimes I$.  If $F$ is orientable, then this $I$-bundle is
simply $F \times I$; in this case we assume an orientation for the base $F$
and the fiber $I$.  We consider the double $F \ltimes S^1$ of $F \ltimes
I$, which again when $F$ is orientable is just $F \times S^1$.

If $p_1,p_2,\dots,p_n$ are points of $F$, then we can apply a Dehn surgery
with slope $b_k/a_k$ to a solid torus neighborhood of the fiber over
$p_k$ in $F \ltimes S^1$, where $a_k$ is a positive integer and $b_k$ is a
relatively prime integer of either sign.  The resulting oriented 3-manifold
$N$ is thus construction from its \emph{Seifert data}, namely the multiset
\[ \{F,(a_1,b_1),(a_2,b_2),\dots,(a_n,b_n)\}. \]
In general we interpret $F$ as an orbifold.  If $a_k \ge 2$, then we
interpret $p_k \in F$ as an orbifold point of order $a_k$, and the circle
over it is an \emph{exceptional fiber}.  By Seifert's classification,
the integers $a_k$ with $a_k \ge 2$ together with the residues $b_k \in
\Z/a_k$ are all topological invariants of the fibration of $N$.  If $F$
and therefore $N$ has boundary, then this is a complete set of invariants.
If $N$ is closed, then the Euler number
\[ e(N) = b + \sum_k \frac{b_k}{a_k} \]
is the only additional necessary invariant.  Thus there is a canonicalized
version of the Seifert data in the form
\[ \{F,b,(a_1,b_1),(a_2,b_2),\dots,(a_n,b_n)\}, \]
where $b$ represents $(1,b)$ and otherwise $a_k \ge 2$ and $1 \le b_k
< a_k$.  If $F$ is non-compact, then $b$ is irrelevant and we omit it
in the canonical form.

With the notation of fiber bundles and Seifert-fibered spaces, we name
these specific manifolds:

\begin{enumerate}
\item We use $S^1 \times S^1$ to denote the standard 2-torus, and $T$
to denote an arbitrary 2-torus, \ie, $T \cong S^1 \times S^1$.
\item $K^2 = S^1 \ltimes S^1$ is the 2-dimensional Klein bottle.
\item $L(m,n)$ is the lens space defined by the Seifert data
$\{S^2,0,(n,m)\}$.
\item $R(m,n)$ denotes the prism space defined by the Seifert data
$\{P^2,0,(m,n)\}$.
\end{enumerate}

\section{Geometrization is recursive}

The goal of this section is to prove \Thm{th:geomrec}, which says that the
geometric decomposition of a 3-manifold $M$ is computable.

\subsection{Statement of geometrization}
\label{ss:geom}

We begin with three results that, together, are one formulation of the
geometrization theorem for closed, oriented 3-manifolds.   

\begin{theorem}[Kneser-Milnor \cite{Kneser:flachen,Milnor:unique}] Every
closed, oriented 3-manifold (other than $S^3$) is a connected sum of prime,
closed, oriented 3-manifolds (none of which are $S^3$).  The summands are
unique up to oriented homeomorphism.
\label{th:factor} \end{theorem}

We will adopt the convenience that a 3-sphere $S^3$ counts as a prime
3-manifold, notwithstanding that \Thm{th:factor} would be easier to state
if $S^3$ were instead interpreted as the ``unit" in the terminology of
unique factorization.

\begin{theorem}[Jaco-Shalen-Johansson \cite{JS:seifert,Johannson:homotopy}] A
closed, oriented, prime 3-manifold has a minimal collection of incompressible
tori, unique up to isotopy and possibly empty, with the property that the
complementary regions are either Seifert-fibered or atoroidal.
\label{th:jsj} \end{theorem}

The decomposition in \Thm{th:jsj} is called the \emph{JSJ decomposition}.
We can call the tori \emph{JSJ tori}, and the complementary regions
\emph{JSJ components}.  We will use $M$ to denote a general closed,
oriented 3-manifold; then $W$ to denote a prime summand of $M$; then $N$
to denote a JSJ component of $W$.

\begin{theorem}[Thurston-Hamilton-Perelman] Suppose that $N$ is an oriented,
prime, atoroidal 3-manifold which is either closed or has torus boundary
components.  Then $N$ is either Seifert-fibered, or it is closed and has a
unique hyperbolic structure, or its interior $N^*$ has a unique, complete
hyperbolic structure with torus cusps.
\label{th:geom} \end{theorem}

As everyone knows, \Thm{th:geom} was conjectured and partly proven by
Thurston \cite{Thurston:kleinian}, then fully proven by Perelman using the
Ricci flow program of Hamilton \cite{MF:ricci}.  (Note that \Thm{th:geom}
implicitly includes the Poincar\'e conjecture in the Seifert-fibered case.)

\begin{remark} Mixing the JSJ decomposition with hyperbolization is a less
pure approach than Thurston's decomposition into geometric components, but
we find it convenient for \Thm{th:main}.   We could recognize spherical
and Euclidean components with the same methods as hyperbolic components
(\Lem{l:hyprec}), while several of the other Thurston geometries induce
canonical Seifert fibrations.  In fact, every Seifert-fibered 3-manifold
or component is geometric.  Conversely, every geometric 3-manifold or
component is hyperbolic unless it is Seifert-fibered or a Sol manifold.
\end{remark}

\subsection{Statement of computational geometrization}
\label{ss:compgeom}

\begin{theorem} If $M$ is a triangulated 3-manifold, then it is recursive
to compute a decorated triangulation which is adapted to its geometric
decomposition.
\label{th:geomrec} \end{theorem}

Before proving \Thm{th:geomrec}, we need to state it more precisely.
When $\Theta$ is a \emph{decorated, adapted triangulation} of $M$ it
means that:
\begin{enumerate}
\item $M$ has a distinguished (but possibly empty) collection of
disjoint, separating 2-spheres, each triangulated with 4 triangles
in $\Theta$, that separates it into prime summands $\{W\}$.  Each $W$
is closed; it inherits its triangulation from $\Theta$ and its holes are
plugged with fresh tetrahedra.
\item The triangulation of each $W$ supports a distinguished (but possibly
empty) collection of disjoint thickened tori $T \times I$ and restricts
to a shelled triangulation of each one.  These thickened tori separate $W$
into JSJ components $\{N\}$.
\item The tetrahedra at all stages are consistently oriented, to express
an orientation of each summand $W$ and each JSJ component $N$ that is
consistent with the orientation of $M$.
\item If $N$ is Seifert-fibered with base $F$, then we make a triangulation
which is adapted to Seifert's description of $N$ by Dehn surgery on $F
\times S^1$ when $F$ is orientable, or Dehn surgery on a twisted
bundle $F \ltimes S^1$ when $F$ is non-orientable.  This includes the case
where $N = W = M$ is a 3-sphere.
\item If $N$ is hyperbolic, then it is marked as the barycentric
subdivision of a regular, adapted cellulation $\Lambda$.  $\Lambda$
comes from a geometric triangulation $\Lambda^*$ of $N^*$ in which each
tetrahedron has at most one ideal vertex.  If $\Delta \in \Lambda^*$ has
an ideal vertex, then it is truncated to a triangular prism in $\Lambda$;
if not, then it is kept in $\Lambda$.  Each tetrahedron in $\Lambda^*$
is also decorated with its dihedral angles.
\end{enumerate}

We proceed to explain each stage of the definition.

\subsubsection{The prime decomposition}

Note that we take a triangulation to be a simplicial complex structure
rather than a generalized triangulation. In geometric topology, for
instance in the SnapPea census, a \emph{generalized triangulation} is
sometimes defined to be a CW-complex whose cells are simplices and whose
attaching maps take simplices to simplices.  In particular, a simplex
in a generalized triangulation need not have distinct vertices and two
simplices may have the same vertices.  We can form a connected sum of two
triangulated 3-manifolds by removing a single tetrahedron from each one
and gluing the sphere boundaries.

\subsubsection{Shelled triangulations}
\label{sss:shelled}

If $X$ is a closed $n$-manifold with two triangulations $\Theta_0$ and
$\Theta_1$, then a \emph{shelled triangulation} of $X \times I$ is a
simplicial complex $\Theta_{0,1}$ whose $(n+1)$-simplices are numbered.
Taken in order, the $(n+1)$-simplices connect the triangulation $\Theta_0$
of $X \times \{0\}$ to the triangulation $\Theta_1$ of $X \times \{1\}$ via
a sequence of bistellar moves.  Note that this combinatorial restriction on
$\Theta_{0,1}$ implies $\Theta_{0,1}$ is PL homeomorphic to $X \times I$.
In other words, if we build $\Theta_{0,1}$ from a sequence of bistellar
moves, and if $X \times \{0\}$ and $X \times \{1\}$ are disjoint in the
result, then $\Theta_{0,1}$ is a triangulation of $X \times I$.

\subsubsection{Orientations}

To be precise, we can decorate each tetrahedron by ordering its vertices,
where two orderings are equivalent if they differ by an even permutation.

\subsubsection{Seifert-fibered components}
\label{sss:seifcomp}

We begin with preliminaries on cellulations and barycentric
subdivisions that we will also need in \Sec{sss:hypcomp}.

A \emph{cellulation} of a topological space $X$ is a CW complex $\Lambda$
with a homeomorphism to $X$.   The complex $\Lambda$ is \emph{regular}
if $\Lambda$ is locally finite; and if the attaching map of each $k$-cell
is an embedding in the $(k-1)$-skeleton, so that each closed $k$-cell of
$\Lambda$ is embedded in $X$.

We will use the following standard proposition to model regular CW complexes
using triangulations.

\begin{proposition}[{\cite[Sec.~10.3.5]{Kozlov:algtop}}] Every regular
CW complex $\Lambda$ has a barycentric subdivision $\Theta$ which is
a simplicial complex, and the spaces of $\Theta$ and $\Lambda$ are
homeomorphic.
\label{p:bary} \end{proposition}

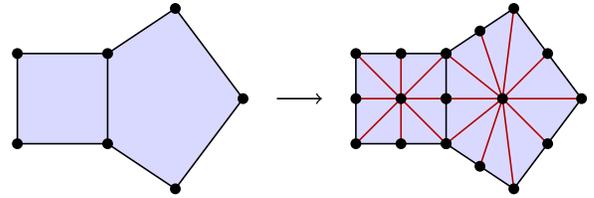
\begin{figure}[htb]\begin{center}\begin{tikzpicture}[scale=.6,semithick]
\draw[fill=lightblue] (1,-1) -- (-1,-1) -- (-1,1) -- (1,1)
    -- (2.5,2) -- (4,0) -- (2.5,-2) -- cycle;
\draw (1,-1) -- (1,1);
\foreach \x/\y in {-1/-1,-1/1,1/1,1/-1,2.5/2,4/0,2.5/-2} {
    \fill (\x,\y) circle (.12); }
\draw[->] (4.75,0) -- (5.75,0);
\begin{scope}[shift={(7.5,0)}]
\draw[fill=lightblue] (1,-1) -- (-1,-1) -- (-1,1) -- (1,1)
    -- (2.5,2) -- (4,0) -- (2.5,-2) -- cycle;
\draw (1,-1) -- (1,1);
\foreach \x/\y in {-1/-1,-1/0,-1/1,0/1,1/1,1/0,1/-1,0/-1} {
    \draw[darkred] (0,0) -- (\x,\y); \fill (\x,\y) circle (.12); }
\fill (0,0) circle (.12);
\foreach \x/\y in {1/1,1.75/1.5,2.5/2,3.25/1,4/0,3.25/-1,2.5/-2,
    1.75/-1.5,1/-1,1/0} {
    \draw[darkred] (2.25,0) -- (\x,\y); \fill (\x,\y) circle (.12); }
\fill (2.25,0) circle (.12);
\end{scope}
\end{tikzpicture}\end{center}
\caption{A barycentric subdivision of part of a cellulation of
    a surface.}
\label{f:bary}\end{figure}

See \Fig{f:bary} for an example.

If $N$ is a Seifert-fibered component, then as described in \Sec{s:notation},
it has a base orbifold $F$ with one circle for each boundary torus of $N$.
The fibration has canonical Seifert data
\[ \{F,b,(a_1,b_1),(a_2,b_2),\dots,(a_n,b_n)\}, \]
with $b$ omitted when $F$ or $N$ has boundary.  The data indicates surgery
with slope $b_k/a_k$ at the fiber over some $p_k \in F$ and (if it exists)
surgery with slope $b$ at $p_0 \in F$.

We choose a triangulation $\Theta_F$ of $F$ such that each $p_k$ (including
$p_0$, if it exists) lies in the interior of a triangle, and such that all
of these triangles are disjoint.  We can lift $\Theta_F$ to a cellulation
$\Lambda$ such that the solid torus $\Delta \times S^1$ over each triangle
$\Delta$ in $\Theta_F$ is tiled by two vertical triangular prisms.  We take
the barycentric subdivision of $\Lambda_F$ to obtain a triangulation of $F
\times I$ or $F \ltimes I$.  If a triangle $\Delta \in \Lambda$ contains
some $p_k$, we remove the solid torus $\Delta \times S^1$ (which is now
triangulated with 72 tetrahedra) and glue it back using Dehn surgery.
The gluing involves a homeomorphism of the boundary $\del (\Delta \times
S^1)$, which we implement with a shelled triangulation of a thickened
torus, as in \Sec{sss:shelled}.  The result is a triangulation of $N$,
which we decorate with information about how it was constructed, so that
the canonical Seifert data is part of the decoration.

\subsubsection{Hyperbolic components}
\label{sss:hypcomp}

If the component $N$ is hyperbolic, then we choose a \emph{geometric
triangulation} $\Theta^*$ of $N^*$, meaning one whose tetrahedra lift
to geometric tetrahedra in the universal cover $\H^3$.  More precisely,
if $\hN$ is the compactification of $N$ given by collapsing each torus
boundary component to a point, we assume a continuous map
\[ f:\Theta^* \to \hN \]
such that the image $f(\Delta)$ of each combinatorial tetrahedron $\Delta$
lifts to a geometric tetrahedron in the standard compactification
$\overline{\H^3}$ of hyperbolic space.  If none of the vertices of
$f(\Delta)$ are at infinity, then $f(\Delta)$ is \emph{finite}; if they
are all at infinity, then $f(\Delta)$ is \emph{ideal}; and if some are at
infinity, then $f(\Delta)$ is \emph{semi-ideal}.   We will assume that all
of the simplices of our $\Theta^*$ are either finite or semi-ideal with one
ideal vertex.  If $N = N^*$ is closed, then all tetrahedra in $\Theta^*$
must be finite; if $N$ has boundary components and thus $N^*$ has cusps,
then some of the tetrahedra must be semi-ideal.

Before proceeding further, we contrast this with some other models that
have also been studied as geometric triangulations.  In some treatments $f$
is not a homeomorphism but only a homotopy equivalence.  In the case we
can still ask for the restriction of $f$ to each tetrahedron $\Delta$ to be
affine-linear in the Klein model of $\overline{\H^3}$.  However, $f(\Delta)$
may be \emph{degenerate}, meaning that it has zero volume, or it may be
\emph{flipped over}, meaning that it has negative signed volume relative
to the orientation of $\Theta^*$ and the standard orientation of $\H^3$.
In another variation, which is often used when $N$ is closed, the inverse
map $g:\hN \to \Theta^*$ is defined, and a finite set of closed geodesic
curves in $\hN$ collapse to ideal points; but the inverse image of any
open tetrahedron in $\Theta^*$ is still a geometric tetrahedron in $\hN$.
Such a structure is a \emph{spun triangulation}, because a geodesic
circle $C \subseteq N$ is approached by cusps of ideal tetrahedra that
wind helically around it.   In particular SnapPea uses spun triangulations.

Ideal geometric triangulations are especially desirable in computational
hyperbolic geometry because they are rigid and algebraically the simplest.
However, it is only a conjecture that every suitable hyperbolic manifold
has an ideal, possibly spun geometric triangulation.  Such a structure
does always exist with degenerate or flipped-over tetrahedra, but these
are less desirable.  We will use finite and semi-ideal tetrahedra in order
to avoid this impasse.  The following proposition is then standard:

\begin{proposition} If $N^*$ is a complete hyperbolic manifold which
is either cusped or closed, then it has a geometric triangulation with
finite and semi-ideal tetrahedra (none of which are spun, degenerate, or
flipped over).  Also, every semi-ideal tetrahedron has only one ideal vertex.
\label{p:triang} \end{proposition}

\begin{proof} We can choose a point $p \in N^*$ and consider the Voronoi
tiling of its orbit in $\H^3$.  Each Voronoi cell is a fundamental domain
and yields a cellulation $\Lambda_1$ of $\hN$.  $\Lambda_1$ is not in
general regular, but it has a barycentric subdivision $\Lambda_2$ which
is regular.  Moreover, each simplex of $\Lambda_2$ has at most one vertex
of $V$ and thus at most one ideal vertex.  We can let $\Theta = \Lambda_3$
be a second barycentric subdivision, which is then a simplicial complex and
still has the property that each tetrahedron has at most one ideal vertex.
\end{proof}

\begin{figure}[htb]\begin{center}\begin{tikzpicture}[thick]
\fill[lightblue] (-2,.5) -- (2,0) -- (.5,3) -- (.125,3.375)
    -- (-.5,3.125) -- cycle;
\draw (-.5,3.125) -- (-2,.5) -- (2,0) -- (.5,3);
\draw[dashed] (-2,.5) -- (.5,1.5) -- (2,0);
\draw[dashed] (.5,1.5) -- (.125,3.375);
\draw[darkred] (-.5,3.125) -- (0,4) -- (.5,3);
\draw[darkred] (.125,3.375) -- (-.5,3.125) -- (.5,3) -- cycle;
\draw[darkred,dashed] (0,4) -- (.125,3.375);
\end{tikzpicture}\end{center}
\caption{A tetrahedron truncated at one vertex.}
\label{f:trunc}\end{figure}
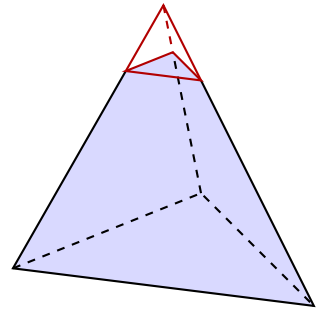

Given a semi-ideal triangulation of $N^*$, we can truncate the cusps so that
each semi-ideal tetrahedron becomes a triangular prism, as in \Fig{f:trunc}.
A barycentric subdivision of this cellulation is then the desired adapted
triangulation.

\subsection{Proof of \Thm{th:geomrec}}

\begin{lemma} It is recursive to find a geometric triangulation of a
hyperbolic 3-manifold $N$ which is either closed or has torus boundary
components, using either of two descriptions of each dihedral angle $0 <
\alpha < \pi$ of each tetrahedron:
\begin{enumerate}
\item Each imaginary exponential $\exp(i\alpha)$ is specified as an element
of $\bQ$.
\item Each angle $\alpha$ is given as a computable real number.
\end{enumerate}
Hence, it is in $\RE$ to determine if $N$ is hyperbolic.
\label{l:hyprec} \end{lemma}

Although the second case of \Lem{l:hyprec} immediately follows from the
first one, we will give a separate proof of each case.  Moreover, even the
weaker second case of \Lem{l:hyprec} is sufficient to prove \Thm{th:main}.

\begin{remark} Manning \cite[Thm.~5.2]{Manning:algo} also proves
\Lem{l:hyprec}, but as a corollary of a harder result.  His results show
(without geometrization) that it is recursive to decide whether $N$ is
hyperbolic, when there is an algorithm for the word problem for $\pi_1(N)$.
He also uses a single polyhedral fundamental domain to describe the
geometry of $N$.  Although this differs from a hyperbolic triangulation,
which is what we use, the two models are somewhat interchangeable for
our purposes.
\end{remark}

\begin{proof}[Proof of case 1 of \Lem{l:hyprec}] Suppose that $\Theta^*$ is a
geometric triangulation of $N^*$.  We can model each tetrahedron $\Delta \in
\Theta^*$ (non-uniquely) by choosing four vertices in the Poincar\'e upper
half-space model, including one on the boundary if $\Delta$ is semi-ideal.
(Note that the ideal vertices of $\Theta^*$ are marked in advance.)
There is an algebraic formula for each finite edge length $\ell$ and each
dihedral angle $\alpha$ of $\Delta$, if these are represented by their
exponential values $\exp(\ell)$ and $\exp(i\alpha)$.  The main matching
condition for $\Theta^*$ to be geometric is that if two tetrahedra share
a finite edge, then the edge lengths agree; and the total dihedral angle
around each edge equals $2\pi$.  The first condition is immediately an
algebraic condition, although note that if an edge is semi-ideal,
then it has infinite length and its length equation is vacuous.

The second condition is almost an algebraic condition since the product of
the exponentiated angles must be 1; this shows that the total angle is a
multiple of $2\pi$, although not which one.  However, this can be remedied
with additional algebraic inequalities, recalling that we are allowed
real algebraic equations for the real and imaginary parts $\cos(\alpha)$
and $\sin(\alpha)$ of each complex variable $\exp(i\alpha)$.  Suppose that
every edge of $\Theta^*$ has at most $n$ incident tetrahedra.  Then we can
make a finite covering of the unit circle $S^1 \subseteq \C$ by rational
rectangles such that each one has an angular extent of less than $2\pi/n$.
We can then loop over choices for which rectangle contains each exponentiated
angle $\exp(i\alpha)$.  If each angle is confined to such a rectangle, we
can know whether the sum of the angles around an edge is specifically $2\pi$
and not some other multiple of $2\pi$.

After forming algebraic equations for all of the geometric data, the
equations have a solution in terms of real algebraic numbers when they
have a solution at all.  For any fixed triangulation $\Theta$, we can thus
use \Thm{th:realalg} (not \Thm{th:tarski}; see the remark
after the proof) to search for a solution and eventually find it, if it
exists.  We must also search over triangulations using \Thm{th:stellar}
or \Thm{th:bistellar}.   Since the result is a nested infinite search
(over triangulations and then candidate geometric structures), we can
apply \Prop{p:gsearch}.
\end{proof}

\begin{remark} Although an infinite search for a solution to algebraic gluing
equations is preposterous in practice, it is good enough for an algorithm
in $\RE$.  Alternatively, for each triangulation, we can apply the more
difficult \Thm{th:tarski} to determine in $\ccR$ if there is a solution.
\end{remark}

\begin{remark} If we allowed geometric triangulations with fully ideal edges,
then it would not be enough for the sum of the angles around such an edge
$e$ to be $2\pi$.  Since $e$ goes to itself under hyperbolic translation
as well as rotation, gluing together the tetrahedra that contain $e$ could
create a non-trivial translational holonomy.  The two conditions together,
that the total angle is $2\pi$ and the translational holonomy vanishes,
are known as a Neumann-Zagier gluing relation \cite{NZ:volumes}.
\end{remark}

\begin{remark} Instead of calculating lengths and angles using positions
of vertices in hyperbolic geometry, we can also relate them directly using
formulas from hyperbolic and spherical trigonometry.
\end{remark}

The separate proof of the second case of \Lem{l:hyprec} works directly
with computable numbers, in effect using numerical analysis to calculate
better and better approximate solutions.  In this approach, we need a
criterion to know that an approximate solution is close to an exact one.
Given a smooth multivariate equation $f(x) = 0$, the Newton-Kantorovich
theorem \cite{Kantorovich:applied} establishes a sufficient criterion for
Newton's method to converge from an approximate solution $x_0$ to an exact
solution $x_\infty$.  Neuberger \cite{Neuberger:monthly} points out that
an ODE analogue of Newton's method, which is called the \emph{continuous
Newton's method}, simplifies the Newton-Kantorovich result.

\begin{theorem}[Newton-Kantorovich-Neuberger
{\cite[Thm.~2]{Neuberger:monthly}}] Let $B_\eps(x_0) \subset \R^n$ be the
open ball of radius $\eps > 0$ around $x_0 \in \R^n$, and let
\[f:B_\eps(x_0) \to \R^n \]
be a $C^2$-smooth function with non-singular matrix derivative $Df$.
Suppose that
\begin{eq}{e:nkn}
||(Df(x))^{-1}f(x_0)|| < \eps
\end{eq}
for all $x \in B_\eps(x_0)$, where $||\cdot||$ is the Euclidean norm
on $\R^n$.  Then there is a unique $x_\infty \in B_\eps(x_0)$ such
that $f(x_\infty) = 0$.  Also, given a solution $x_\infty$ such that
$Df(x_\infty)$ is non-singular, equation \eqref{e:nkn} eventually
holds as $x_0 \to x_\infty$, moreover with $\eps \to 0$.
\label{th:nkn} \end{theorem}

Although we will not reprove \Thm{th:nkn}, we can discuss where the theorem
and its proof come from.  Newton's method to find a root of a univariate
function $f:(a,b) \to \R$ begins at an approximate root $x_0 \in (a,b)$
and applies the iteration
\[ x_{n+1} = x_n - \frac{f(x_n)}{f'(x_n)}, \]
which in favorable cases converges to a solution $x_\infty$ of the equation
$f(x) = 0$.  If $f$ is multivariate as in \Thm{th:nkn}, then this has the
well-known matrix generalization
\[ x_{n+1} = x_n - (Df(x_n))^{-1}f(x). \]
Finally in the continuous version, we let $x(0) = x_0$ and define
the ODE
\[ x'(t) = -(Df(x(t)))^{-1}f(x). \]
Then in favorable cases the limit
\[ x_\infty = \lim_{t \to \infty} x(t) \]
is again a solution to $f(x) = 0$.

\begin{remark} Although Neuberger's paper on the continuous Newton's method
is more recent than Thurston's work, Kantorovich's earlier, more complicated
formula also suffices for \Lem{l:hyprec} and \Thm{th:main}.
\end{remark}

If the equation $f(x) = 0$ has a non-singular Jacobian $Df$ in a neighborhood
of a solution, as in \Thm{th:nkn}, then the system of equations is also
called \emph{transverse} or \emph{first-order rigid}.  We will need
a generalization of this concept.  Given a smooth function
\[ U \subseteq \R^n \qquad f:U \to \R^m, \] 
where $n$ and $m$ need not be equal, if $Df$ has constant rank $kr$, then the
image $f(U)$ is a manifold and $f$ is a submersion onto its image.  In this
case the equation $f(x) = 0$ is first-order rigid \emph{except for} the
directions parallel to the manifold $f^{-1}(0)$.  By the implicit function
theorem, we can discard some set of $n-rk$ coordinates in the domain and
project to some set of $k$ coordinates in the target to achieve unconditional
first-order rigidity that satisfies the hypotheses of \Thm{th:nkn}.

To establish first-order rigidity in our case, we will need a corollary
of the Calabi-Weil rigidity theorem.

\begin{theorem}[Calabi-Weil {\cite[Sec.~8.10]{Kapovich:book}}]  If $N$
is a closed, hyperbolic 3-manifold, then the induced representation of
its fundamental group,
\[\rho:\pi_1(N) \to \Isom(\H^3), \]
is first-order rigid except for conjugacy.  (\Ie, it is infinitesimally
rigid at the level of the first derivative.)  The same is true if $N$
is cusped, among representations that are parabolic at the torus cusps.
\label{th:rigid} \end{theorem}

\begin{corollary}[Stated by Izmestiev {\cite[Sec.~1.3]{Izmestiev:rigidity}}]
If $\Theta$ is a geometric triangulation of a closed or cusped hyperbolic
3-manifold $N^*$, then it is first-order rigid except for motion of the
non-ideal vertices.
\label{c:rigid} \end{corollary}

Since we could not find a proof of \Cor{c:rigid} in the literature, we
provide one in \Sec{ss:rigid}.

\begin{proof}[Proof of case 2 of \Lem{l:hyprec}] We fix the model of
each tetrahedron in the upper half space model so that it has exactly six
degrees of freedom, or five if one of the vertices is ideal.  After ordering
the vertices $v_1, v_2, v_3, v_4$, we can put vertex $v_1$ at $(0,0,1)$;
we can put vertex $v_2$ directly below it (or at $(0,0,0)$, allowing it
to be the ideal vertex); and we can put $v_3$ at a position of the form
$(a,0,b)$.  We approximate the positions of the vertices approximately
with rational numbers.  We can then approximate the lengths and angles
of each tetrahedron in the same form, as well as the first and second
derivatives of the lengths and angles as a function of the main variables,
the separate positions of the vertices in the ideal models of the tetrahedra.

Suppose that there are $n$ non-ideal vertices.   By the implicit function
theorem, any exact solution to the gluing equations for the tetrahedra can
be perturbed so that some $3n$ of the coordinates are exactly rational.
Also by the implicit function theorem, some $3n$ of the angle conditions
are implied by the other angle conditions and can be omitted.  Finally,
the fixed coordinates and omitted angle conditions can be chosen so that
the remaining system of constraints, which we can write abstractly as $f(x)
= 0$, has a non-singular Jacobian $Df$.

Moreover, the mapping $f$ is real analytic with an explicit formula.  Thus,
given an approximate solution $x_0$ which is within $\eps$ of a true solution
and $\eps$ is small enough, we can majorize $||(Df(x))^{-1}||$ on the ball
$B_\eps(x_0)$ using Taylor series, to confirm equation \eqref{e:nkn}.
\end{proof}

As \Lem{l:hyprec} addresses the hyperbolic case of \Thm{th:geomrec},
we turn to a lemma and a proposition that address the Seifert-fibered case.

\begin{lemma} It is recursive to find an adapted triangulation of a
Seifert-fibered manifold $N$ which is either closed or has torus boundary
components.  Hence, it is in $\RE$ to determine if it is Seifert-fibered.
\label{l:seifrec} \end{lemma}

\begin{proof} We can search through triangulations until we find one
that is a barycentric subdivision of a cellulation by triangular prisms.
It is then easy to check whether the prisms fit together following the
rules in \Sec{ss:compgeom}.
\end{proof}

Finally, a torus $T$ that has matching Seifert-fibered structure on both
sides is not needed and is not a JSJ torus.   It is easy to see
this case in the proof of \Thm{th:jsj}.  The more subtle possibility
is that one or both sides might have more than one Seifert fibration.
Fortunately this is rare for Seifert-fibered manifolds with boundary.
It is addressed by the following result.

\begin{theorem}[Waldhausen
{\cite[Thm.~VI.17 \& Lem.~VI.19]{Jaco:lectures}}] Let $N$ be an oriented
3-manifold with non-empty boundary $\del N$ and which has at least
one Seifert fibration.  Then the fibration is uniquely determined up to
isotopy by its restriction to $\del N$, and is outright unique except the
following cases:
\begin{enumerate}
\item If $N$ is a solid torus $D^2 \times S^1$, then every fibration
of $\del N$ extends to a fibration of $N$ over a disk $D^2$
with at most one exceptional fiber.
\item If $N$ is a thickened torus $S^1 \times S^1 \times I$, then every
fibration is a trivial circle bundle over an annulus.  There is such a
fibration for every rational slope in a single torus $S^1 \times S^1$.
\item If $N$ is a twisted $I$-bundle $K^2 \ltimes I$ over a Klein bottle
$K^2$, then it has two non-isotopic fibrations.  One fibration is over a
M\"obius strip with Seifert data $\{S^1 \ltimes I\}$, and one is over a
disk $D^2$ with Seifert data $\{D^2,(2,1),(2,1)\}$.
\end{enumerate}
\label{th:waldhausenb} \end{theorem}

\begin{proof}[Proof of \Thm{th:geomrec}] We search over triangulations
$\Theta$ of $M$ using stellar or bistellar moves, and decorations of them, to
find an adapted triangulation as described in \Sec{ss:compgeom}.  A suitable
decoration consists of distinguished spheres and thickened tori, and a
reverse barycentric subdivision in each JSJ component $N$ to make triangular
prisms in the Seifert-fibered case and a combination of once-truncated
and ordinary tetrahedra in the hyperbolic case.  Within this search, we
search for geometric data to describe the hyperbolic structure of each $N$
which is not Seifert-fibered.  Since these are nested, infinite searches,
we combine them using the $\RE$ search algorithm of \Prop{p:gsearch}.
By the geometrization theorem, we will eventually find a $\Theta$ that
fits the description of \Sec{ss:compgeom}.

We examine the JSJ components to verify that all of the spheres and tori are
essential and that no two are parallel.  We veto $\Theta$ if it includes
a Seifert-fibered solid torus (by checking that the base orbifold is a
disk with at most one exceptional fiber).  If all JSJ components are either
hyperbolic or Seifert-fibered, and if none are solid tori, then every sphere
and torus is essential and no two spheres are parallel.  Two distinct tori
are parallel if and only if the component between them is a thickened torus;
we veto this as well.  We also need to veto a Seifert-fibered component
homeomorphic to $P^3 \sharp P^3$, which is the only Seifert-fibered space
that is not a prime 3-manifold.  This space has Seifert data $\{P^2,0\}$.

Finally, we need to check that all of the tori are JSJ tori.  We veto
$\Theta$ if there is a torus $T$ that has Seifert-fibered components
on both sides that:
\begin{enumerate}
\item restrict to the same fibration of $T$; or
\item could be refibered to restrict to the same fibration of $T$.
\end{enumerate}
Case 1 is easy to recognize.   By \Thm{th:waldhausenb} and the comments
after, in case 2 we only have to consider two types of Seifert-fibered
components, which can be recognized explicitly from any of their fibrations:
\begin{enumerate}
\item[2a.] $N \cong S^1 \times S^1 \times I$, or 
\item[2b.] $N \cong K^2 \ltimes I$.
\end{enumerate}
Case 2a is only possible if $N$ is glued to itself to make $W$ a torus bundle
over a circle, $S^1 \ltimes (S^1 \times S^1)$, because we do not allow
parallel tori.  The torus is needed if and only if $W$ is a Sol manifold.
We can verify this case by confirming that the holonomy matrix in $\SL(2,\Z)$
has distinct, real eigenvalues.

In case 2b, $N \cong K^2 \ltimes I$ only has one torus boundary component
$T$, so its refibration does not affect any other torus.  In this
case the fibration of $N$ may have Seifert data $\{S^1 \ltimes I\}$ or
$\{D^2,(2,1),(2,1)\}$.  The resulting binary choice may occur on one or
both sides of $T$, and we veto $\Theta$ if the Seifert fibrations extend
across $T$ for any of these choices.
\end{proof}

\subsection{Proof of \Cor{c:rigid}}
\label{ss:rigid}

The idea of the proof is that we can convert a first-order deformation
of a triangulation of $N^*$ to a deformation of a representation of $\rho$,
in much the same way that we can convert a triangulation to $\rho$ in the
first place.

\begin{proof}
In general, if $\Gamma$ is a discrete group (such as the fundamental group
of a topological space) and $G$ is a Lie group, then we can describe
a first-order deformation of a homomorphism $\rho:\Gamma \to G$ as a
homomorphism
\[ (\rho,\rho'):\Gamma \to G \ltimes \mg. \]
Here $\mg$ is the Lie algebra of $G$ viewed as a group under addition,
while $G \ltimes \mg$ is the semidirect product in which the non-normal
subgroup $G$ acts on the normal subgroup $\mg$ by conjugation.
Also, $(\rho,\rho')$ should reduce to $\rho$ under the quotient map
\[\pi:G \ltimes \mg \to G. \]
Note that $G \ltimes \mg$ is also the total space of the tangent bundle
$TG$.  In other words, the extension $\rho'$ is a choice of a tangent
vector $\rho'(g) \in T_{\rho(g)}G$ for every $g \in \Gamma$, such that
the pairs $(\rho(g),\rho'(g))$ together make a group homomorphism.

Suppose that $\Gamma = \Gamma_1(X)$ is the fundamental group of a based CW
complex $X$.  Then we can model $\rho$ (non-uniquely) as a non-commutative
cellular cocycle $\alpha \in C^1(X;G)$.  Given $\rho$, we can likewise
model the extension $\rho'$ (also non-uniquely) as a commutative cocycle
$\alpha' \in C^1(X;\mg)$, where here $\mg$ is a coefficient system twisted
by $\alpha$.

Now let $X = N$, where $N$ has a cellulation $\Theta$ that comes from a
closed or cusped hyperbolic structure on $N^*$ and a geometric triangulation
$\Theta^*$.  If $N$ is cusped and $\Theta^*$ is a semi-ideal triangulation,
then we make $\Theta$ by truncating the ideal vertices of $\Theta^*$.
We then want to make a cocycle $\alpha$ from $\gamma$.   To do this,
we first choose a specific isometry $\tN^* \cong \H^3$.  Then we choose
an orthonormal tangent frame at each vertex of $\Theta$.  Given an edge
$e \in \Theta$, we let $\alpha(e)$ be the element of $G = \Isom^+(\H^3)$
that takes the tail $\tv$ of a lift $\te$ to the head $\tw$, and takes the
lifted frame of $\tv$ to the lifted frame of $\tw$.  If $N^*$ is cusped,
then we require that each truncation edge in $N$ is assigned a parabolic
element that fixes the corresponding ideal vertex in $N^*$.

In this setting, \Thm{th:rigid} says that $H^1(N;\mg) = 0$ in the closed
case and $H^1(N,\partial N;\mg,\mp) = 0$ in the cusped case, where $\mg$
is the parabolic Lie subalgebra of $\mg$.  The theorem is typically proved
using de Rham cohomology rather than cellular cohomology, but these models
of cohomology are isomorphic as usual.  More explicitly, every 1-cocycle
$\alpha' = \delta\beta$, where $\beta$ is an $\mg$-valued 0-cochain on
the vertices of $\Theta$.

Finally, suppose that $\gamma'$ is a first-order deformation of the
hyperbolic structure $\gamma$ of $\Theta^*$ that satisfies the first
derivative of the gluing equations.  Then we can lift $\gamma'$ to a cocycle
$\alpha'$ (non-uniquely) in the same way that $\gamma$ lifts to $\alpha$.
Then \Thm{th:rigid} provides $\beta$, and $\beta$ descends to a first-order
motion of the vertices of $\Theta^*$ that induces the deformation $\gamma'$.
\end{proof}

\section{Homeomorphism is recursive}
\label{s:homrec}

The goal of this section is to prove \Thm{th:main} in this section,
postponing only the proof of \Thm{th:hyphom} below until \Sec{s:hyphom}.

\subsection{Connected sums}
\label{ss:sums}

If $M_1$ and $M_2$ are two closed, oriented 3-manifolds given by
triangulations, then by \Thm{th:geomrec}, we know the direct sum
decompositions of each one into prime 3-manifolds.  These summands can
be freely permuted and can only be matched in finitely many ways.  If we
search over the ways to match them, we then reduce the oriented homeomorphism
problem $M_1 \congq M_2$ to the oriented homeomorphism problem $W_1 \congq
W_2$ for prime summands $W_1$ and $W_2$.  To review, each summand $W_k$
inherits an orientation from its parent $M_k$; in the reverse direction,
there is no ambiguity in forming an oriented connected sum.

\subsection{One JSJ component}

We switch to the other end of geometric decompositions to analyze a single
pair of JSJ components $N_1 \subseteq W_1 \subseteq M_1$ and $N_2 \subseteq
W_2 \subseteq M_2$.  We are interested not only in the isomorphism problem,
but also in the effect of the mapping class group of a component $N$ 
on the boundary $\del N$.

\begin{theorem} Suppose that $N$ is an oriented, hyperbolic JSJ summand
such that $N^*$ is either closed or cusped.  Then the mapping class group
of $N$ is its isometry group.  It is a finite group and its computation
is recursive.  If $N_1$ and $N_2$ are two such manifolds, then they are
homeomorphic if and only if they are isometric, and recognizing this
condition is recursive.
\label{th:hyphom} \end{theorem}

Again, we will prove \Thm{th:hyphom} in \Sec{s:hyphom}.  Note that if $N$
is hyperbolic and has torus boundary components, then each such component
inherits a Euclidean structure from the hyperbolic structure on $N^*$.

Suppose instead that $N$ is Seifert-fibered (and, as usual, oriented).
Then in the direct sense the automorphism problem only matters for
\Thm{th:main} when the JSJ graph is non-trivial and thus $N$ has boundary.
However, we will learn the relevant automorphism properties from an
associated closed Seifert-fibered space.

\begin{lemma} Let $N$ be a closed, oriented 3-manifold which is decorated
with a Seifert fibration with Seifert data
\[ \{F,b,(a_1,b_1),(a_2,b_2),\dots,(a_n,b_n)\}. \]
Then:
\begin{enumerate}
\item The exceptional fibers of $N$ are freely permutable by automorphisms
of the Seifert fibration, provided that the permutation preserves the
orbifold number $a_k \ge 2$ and the residue $b_k \in \Z/a_k$ of each
exceptional fiber.
\item Any finite set of regular fibers is freely permutable.
\item If the base $F$ is orientable, then $N$ has a homeomorphism that
inverts all fibers together, but they cannot be inverted separately.
\item If the base $F$ is non-orientable, then given two disjoint finite
sets $A, B \subseteq F$, $N$ has a homeomorphism that inverts the fibers
over $A$ in place and fixes the fibers over $B$.
\end{enumerate}
\label{l:seifperm} \end{lemma}

\begin{proof} Cases 1, 2, and 4 can all be established by isotopies of $F$
that move points that correspond to the distinguished fibers.  In case 4,
using the hypothesis that $F$ is non-orientable, we can move a point $p
\in A$ around an orientation-reversing loop in $F$ that stays away from $B$
and from the rest of $A$.

Meanwhile in case 3, the fibration itself is orientable, which means that an
orientation of any one fiber induces a canonical orientation of all fibers.
On the other hand, Seifert's construction of the fibration via vertical Dehn
surgery on $F \times S^1$ is invariant with respect to inverting the $S^1$
factor and simultaneously applying an orientation-reversing homeomorphism
to $F$.
\end{proof}

We also need the counterpart to \Thm{th:waldhausenb} for closed
Seifert-fibered spaces.

\begin{theorem}[Waldhausen {\cite[Thm.~VI.17]{Jaco:lectures}}]
If $N$ is a closed, oriented Seifert-fibered 3-manifold, then its Seifert
fibration is unique up to homeomorphism except in the following cases:
\begin{enumerate}
\item A Seifert-fibered space with base $S^2$ and at most two exceptional
fibers is either a lens space $L(m,n)$, $S^2 \times S^1$, or $S^3$
\item A Seifert-fibered space with base $P^2$ and at most one exceptional
fiber is either a lens space $L(4,n)$, a prism space $R(m,n)$, or $P^3
\sharp P^3$.
\item The space with Seifert data
\[ \{S^2,b,(2,1),(2,1),(a_1,b_1)\} \]
is a prism space $R(m,n)$.
\item The twisted bundle $K^2 \ltimes S^1$ which is the double
of $K^2 \ltimes I$ has the double of its two fibrations, namely
the Seifert data $\{K^2,0\}$ and the Seifert data
\[ \{S^2,0,(2,1),(2,1),(2,1),(2,1)\}. \]
\end{enumerate}
\label{th:waldhausenc} \end{theorem}

\Thm{th:waldhausenc} comes with simple formulas for which lens space or
prism space is obtained, which we omit.  In particular, the answer is
recursive and (as we will later want) elementary recursive.

\subsection{The JSJ graph}
\label{ss:graph}

If $W$ is a prime 3-manifold, then its JSJ decomposition is modelled
by a labelled graph $\Gamma$, whose vertices represent JSJ components
and whose edges represent connecting tori.  Each vertex is labelled by
the homeomorphism type of its component, which is either Seifert-fibered
or hyperbolic.  In addition, each edge is decorated with gluing data
and peripheral data which will be described precisely in the proof
of \Thm{th:main} below.

\begin{remark} This graph structure inspired the term \emph{graph manifold}
for a prime 3-manifold whose JSJ components are all Seifert-fibered
\cite{Waldhausen:gruppen}.   This terminology is standard but ironic,
since geometrization shows that the same graph concept is important for
all prime 3-manifolds.
\end{remark}

The labelled graph $\Gamma$ is an invariant of $W$, which at first glance may
seem like a complete invariant, provided that the homeomorphism problem for
each JSJ component is recursive.  However, it is not that simple, because we
have to know the allowed permutations of the torus boundary components
of a JSJ component $N$, and the allowed homeomorphisms of each torus
boundary component.  Finally, we need to deal with the special case that
$N$ is either $K^2 \ltimes I$ or $S^1 \times S^1 \times I$ and thus has
more than one Seifert fibration.

\begin{proof}[Proof of \Thm{th:main}] (Proof using case 1 of \Lem{l:hyprec}.)
As explained in \Sec{ss:sums}, it suffices to solve the homeomorphism
problem $W_1 \congq W_2$ for prime 3-manifolds $W_1$ and $W_2$.  The proof is
divided into three steps.  In steps 1 and 2, we let $W$ be a prime 3-manifold
and let $\Gamma$ be its JSJ graph.  It is recursive to calculate $\Gamma$
and the isomorphism types of its vertices.

Step 1: We address the cases in which a JSJ component of $W$ has more than
one Seifert fibration. If any component is a $K^2 \ltimes I$, then its two
fibrations (described in \Thm{th:waldhausenb}) are inequivalent; we choose
one of them and use it for every occurrence of $K^2 \times I$ in $W$.
If a JSJ component $N$ is a thickened torus, then as in the proof of
\Thm{th:geomrec}, $W$ is a Sol manifold and a torus bundle over a circle,
$S^1 \ltimes (S^1 \times S^1)$.  In this case the homeomorphism type of
$W$ is given by a pair of conjugacy classes in $\SL(2,\Z)$, one for each
orientation of the base circle.  Recall that the conjugacy classes in
$\SL(2,\Z)$ can be classified with the aid of the isomorphism
\[ \PSL(2,\Z) \cong C_2 * C_3. \]
If $g \in \SL(2,\Z)$ has non-zero trace (which it does if $W$ is Sol),
then its conjugacy class is given by the sign of its trace and its reduced
cyclic word in $C_2 * C_3$.

Step 2: We suppose that $W$ is not a Sol torus bundle over a circle.  If $T$
is a JSJ torus in $W$ and one side of $T$ is a hyperbolic component $N$,
then $T$ inherits a Euclidean structure which we can normalize to have
area 1.   This Euclidean structure can be described by a quadratic form $Q$
on the first homology $H_1(T) = H_1(T;\Z)$, where $Q(c)$ is the square of
the minimum length of $c \in H_1(T)$.  Moreover, the coefficients of $Q$
are real algebraic numbers computable from the geometry of $N$.  On the
other hand, if $N$ is Seifert-fibered, then the induced fibration of $T$
selects a \emph{line} in $H_1(T)$, by which we mean a rank-one subgroup
$L \subseteq H_1(T)$ with a torsion-free quotient $H_1(T)/L$.

Since $T$ has two sides, it is then decorated by a pair of quadratic forms
on $H_1(T)$, or a quadratic form and a line, or a pair of lines.  In the
third case when both sides of $T$ are Seifert-fibered, the fibrations
must be mismatched at $T$, so the two lines $L_1, L_2 \subseteq H_1(T)$
must be distinct.  Hence they constitute a \emph{rational line basis}
in the sense that
\[ H_1(T;\Q) = (L_1 \tensor \Q) \oplus (L_2 \tensor \Q). \]
We also obtain a rational line basis in the second case, when one side
is Seifert-fibered and produces a line $L_1 = L$, and the other side is
hyperbolic and produces a quadratic form $Q$.   In this case, there exist
a finite set of pairs of homology classes $\pm c \in H_1(T) \setminus L_1$
that minimize $Q(c)$.  If there is only one such pair, we let $L_2$ be the
line generated by $\pm c$.  If there is more then one, we let $L_2$ be the
line generated by the first such pair in the clockwise direction from $L_1$.

Note that each possible decoration of $T$ induced by the geometry on both
of its sides has a finite stabilizer in the oriented mapping class group
$\SL(H_1(T)) \cong \SL(2,\Z)$ of $T$.  If we order the two sides of $T$,
then the stabilizer usually has two elements; in rare cases it is a cyclic
group of order 4 or 6.

If $N$ is a Seifert-fibered component of $W$, then each of its torus
boundary components is decorated by a rational line basis.  We can thus
make a closed Seifert-fibered space $\hN$ by collapsing a circle fibration
of each component $T \subseteq \del N$ that represents the opposite line
in $H_1(T)$, the one that does not come from $N$ itself.  Each
torus component of $\del N$ becomes a distinguished fiber in $\hN$
which may be either regular or exceptional.

Step 3:  Suppose that $W_1$ and $W_2$ are two prime, closed, oriented
3-manifolds.  If they do not have any JSJ tori, then they are both
Seifert-fibered, and we can use \Thm{th:waldhausenc} to tell if they are
the same.  If they are both Sol torus bundles, we can use the algorithm
in step 1 to determine if they are homeomorphic.

Otherwise we can assume that $W_1$ and $W_2$ have non-trivial JSJ graphs
$\Gamma_1$ and $\Gamma_2$, and that each JSJ torus has a canonical decoration
as described in step 2.  To determine if $W_1$ and $W_2$ are homeomorphic,
we search over graph isomorphisms $f:\Gamma_1 \to \Gamma_2$.  For every pair
of $T_1 \subseteq W_1$ and $T_2 \subseteq W_2$ that are matched by $f$,
we search over mapping classes that preserve the canonical decorations of
$T_1$ and $T_2$.   In the innermost part of the search, we want to calculate
whether the homeomorphisms of the JSJ tori extends to each matched pair
of JSJ components $N_1 \subseteq W_1$ and $N_2 \subseteq W_2$.  If $N_1$
and $N_2$ are hyperbolic, then we can use \Thm{th:hyphom} to determine if
the homeomorphism $\del N_1 \cong \del N_2$ extends.  If they are both
Seifert-fibered, then we can use \Lem{l:seifperm} to determine whether
the corresponding closed Seifert-fibered manifolds $\hN_1$ and $\hN_2$
have a homeomorphism that extends the given homeomorphism $\del N_1 \cong
\del N_2$.  Note that we can employ \Lem{l:seifperm} because any relevant
homeomorphism $N_1 \cong N_2$ preserves the fibration at the boundary, and is
thus isotopic to a fibration-preserving homeomorphism by \Thm{th:waldhausenb}.
\end{proof}

\begin{proof}[Proof of \Thm{th:main}] (Proof using case 2 of \Lem{l:hyprec}.)
If the geometric data of each hyperbolic JSJ component $N$ of a summand
$W$ is described with computable real numbers rather than real algebraic
numbers, then the induced Euclidean structure on a JSJ torus $T \subseteq
\del N$ is only given by a convergent sequence of approximations.  Thus,
it is not possible to definitively calculate the isometries of $T$ or the
shortest cycles, as expressed with the quadratic form $Q(c)$.  However, all
non-isometries and all non-zero classes in $H_1(T)$ that are not shortest
are eventually revealed.  This yields an algorithm in $\coRE$ for the
homeomorphism problem $M_1 \congq M_2$, which is enough to show that the
problem is recursive per the discussion at the beginning of \Sec{s:hyphom}.
\end{proof}

\section{Proofs of \Thm{th:hyphom}}
\label{s:hyphom}

In this section we will give several proofs of \Thm{th:hyphom}.  Recall that
\Cor{c:plre} says that the existence of a PL homeomorphism $N_1 \cong N_2$
is in $\RE$; it is also easy to check whether it preserves orientation.  So,
by \Prop{p:recore}, it suffices to show that homeomorphism is in $\coRE$,
although only one of the proofs will make use of this directly.  By a similar
argument, finding elements in the mapping class group of a single $N$ is in
$\RE$; the remaining task is an algorithm to show that the list is complete.

Recall that if $N$ has boundary, then its interior $N^*$ is cusped and
has a semi-ideal triangulation $\Theta^*$.   In this case, $\Theta$
is a cellulation in which semi-ideal tetrahedra are once truncated.
We want to geometrize the truncation that produces $\Theta$.  We consider
a horospheric truncation which is almost but not quite unique, with the
following three properties:
\begin{enumerate}
\item The horosphere sections lie entirely within the semi-ideal
tetrahedra of $\Theta^*$, and therefore do not intersect each other.
\item For some common integer $n$, every horospheric
torus has area $2^{-n}$.
\item We do not use the smallest value of $n$ that satisfies
conditions 1 and 2.
\end{enumerate}
For convenience, we let $N^* = N$ and $\Theta^* = \Theta$ if $N$ is closed.

Some of the proofs make use of the following lemma.

\begin{lemma} It is recursive to obtain a lower bound
in the injectivity radius of $N$ and $\Theta$.
\label{l:inject} \end{lemma}

\begin{proof}[First proof] Suppose first that $N$ is closed.  For each vertex
$v \in \Theta$, let $U_v$ be the open star of $\Theta$ containing $p$.
Then the collection $\{U_v\}$ is a finite open cover of $N$.  It follows
just from topology that there is some radius $\eps$ such that every ball
of radius $\eps$ is contained in some $U_v$.  For an explicit calculation,
let $\Theta'$ be a barycentric subdivision of $\Theta$, and for each $v \in
\Theta$, let $X_v$ be the closed star of $v \in \Theta'$; then the sets $X_v$
are a closed cover.   We can calculate or bound the distance from $X_v$
to $N \setminus U_w$ for some $w \in \Theta$ with $X_v \subseteq U_w$.
The minimum of all of these distances is thus a lower bound $\eps$ for
the injectivity radius.
\end{proof}

\begin{proof}[Second proof] In general we use the notation $B(p,r)$ for
a hyperbolic ball of radius $r$ centered at $p$.

Let $r$ be the exact injectivity radius of $N$, and let $p$ be a point
on a closed geodesic of $N$ of length $2r$.  Then $p \in \Delta$ for some
cell $\Delta \in \Theta$, and we can let $\ell$ be an upper bound of the
diameter of $\Delta$.  Then in the universal cover
\[ \tN \subseteq \tN^* \cong \H^3, \]
we obtain that at least $1/{2r}$ lifts of $\Delta$ intersect $B(p,1)$, and
thus at least this many copies of $\Delta$ are contained in $B(p,\ell+1)$.
Thus
\[ \frac{1}{2r} \le \frac{\Vol(B(p,\ell+1))}{\Vol(\Delta)}, \]
hence
\begin{eq}{e:ibound}
r > \frac{\Vol(\Delta)}{2\Vol(B(p,\ell+1))}.
\end{eq}
We can calculate an upper bound of this form, if necessary using a lower
bound for the numerator and an upper bound for the denominator, for every
cell in $\Theta$, since we do not know the position of the shortest geodesic
loop in advance.
\end{proof}

\begin{proof}[Third proof] This proof is a variation of the second
proof using the entire diameter and volume of $N$.  J{\o}rgensen and
Thurston proved that the set of possible volumes of $N^*$ is well-ordered.
In particular, there is one of least volume, so there is some constant $c >
0$ such that
\[ \Vol(N^*) > c. \]
Our construction of the geometry of $N$ spares more than half of the volume
of $N^*$, so
\[ \Vol(N) > \frac{\Vol(N^*)}2 > \frac{c}2 = c'. \]
We can obtain an upper bound $\ell$ on the diameter of all of $N$ by
adding bounds on the diameters of the cells in $\Theta$.  Then, we let $D$
be a convex fundamental domain for $N$; it has the same volume and diameter
at most $2\ell$.  Thus we obtain an estimate similar to \eqref{e:ibound},
but more robust:
\[ r > \frac{\Vol(D)}{2\Vol(B(p,\ell+1))}
    > \frac{c'}{2\Vol(B(p,\ell+1))}. \qedhere \]
\end{proof}

\begin{remark} Without an explicit bound on least-volume closed or
cusped hyperbolic manifold, the third proof has the unusual feature of
non-constructively proving that an algorithm exists, \ie, without fully
stating the algorithm.   Meyerhoff \cite{Meyerhoff:lower} established the
first lower bound
\[ \Vol(N) \ge \frac{2}{5^4} \]
in the closed case.  In the same paper, he and J{\o}rgensen established
\[ \Vol(N^*) \ge \frac{\sqrt{3}}4 \implies \Vol(N) \ge \frac{\sqrt{3}}8 \]
in the cusped case.  The exact minimum values are now known
\cite{GMM:minimum}.
\end{remark}

\begin{proof}[First proof of \Thm{th:hyphom}] This proof is similar to one
given by Scott and Short \cite{SS:problem}.  We assume geometric
triangulations $\Theta^*_1$ and $\Theta^*_2$ of $N_1^*$ and $N_2^*$.

If $N^*_1$ and $N^*_2$ are homeomorphic and therefore isometric, then we
can intersect the tetrahedra of $\Theta^*_1$ and $\Theta^*_2$ to make a
tiling of $N_1 \cong N_2$ by various convex cells with 8 or fewer sides;
we can then take a barycentric subdivision to make tetrahedra.   We thus
obtain a mutual refinement $\Theta_3$ of $\Theta_1$ and $\Theta_2$.
If we can bound the complexity of $\Theta_3$, then we can find it with a
finite search or show that it does not exist, rather than using stellar
or bistellar moves in both the up and down directions.

Let $\Delta_1 \in \Theta^*_1$ and $\Delta_2 \in \Theta^*_2$ be two tetrahedra
in the separate triangulations.  In the universal cover $\tN^*_1$, they can
only intersect in a single cell with at most 8 sides.  In $N^*_1$ itself
they can intersect many times; however, only as often as different lifts
of $\Delta_1$ intersect one fixed lift of $\Delta_2$.  If $\Delta_1$ and/or
$\Delta_2$ are semi-ideal, then their lifts intersect if and only if their
truncations do.  There is a recursive volume bound on the number of possible
intersections by the same argument as the second proof of \Lem{l:inject}.

Having bounded the necessary complexity of a mutual refinement $\Theta_3$,
we can now search over separate refinements $\Theta_3$ of $\Theta_1$
and $\Theta_4$ of $\Theta_2$ using \Prop{p:refine}, and look for an
orientating-preserving simplicial isomorphism $\Theta_3 \cong \Theta_4$.
The same method can be used to calculate the mapping class group of a
single $N$.
\end{proof}

\begin{proof}[Second proof] Suppose that $X_1$ and $X_2$ are two compact
metric spaces, and suppose that for each $\eps > 0$ we have a way to make
finite $\eps$-nets $S_1$ and $S_2$ for $X_1$ and $X_2$, and calculate
or approximate all distances within $S_1$ and within $S_2$.  If $X_1$
and $X_2$ are isometric, then there is a function $f:S_1 \to S_2$ that
changes distances by at most $2\eps$.  On the other hand, if there is such
a function for every $\eps$, then $X_1$ and $X_2$ must be isometric.

In our case, we let $X_k = N_k$, where we make sure to use the same
truncation area $2^{-n}$ to geometrize $N_1$ and $N_2$ given the geometries
of $N^*_1$ and $N^*_2$.  We calculate a common lower bound $\delta$ on
the injectivity radius.

We can choose some convenient coordinates inside each cell $\Delta \in
\Theta_k$.  We then have the ability to calculate geodesic segments in $N_k$
that are made of geodesic segments in the separate tetrahedra.  If $\Delta$
is truncated, then the geodesic segment might hug the truncation boundary
for part of its length, but it still has a finite description.  Without
more work, we don't know which of these geodesics are shortest geodesics.
However, if a geodesic is shorter than $\delta$, then it is shortest.
Taking $\delta \gg \eps \to 0$, we can make $\eps$-nets of both $N_1$
and $N_2$ and look for approximate isometries between these $\eps$-nets;
it suffices to check distances below the fixed value $\delta$.

More explicitly, we can use the covering by open stars $S_v$ in the
first proof of \Lem{l:inject}.   There is a $\delta$ such that
if $d(x,y) < \delta$, then $x$ and $y$ and even the connecting short
geodesic are all in some open star.

This algorithm does not by itself ever prove that $N_1$ and $N_2$ are
isometric, only that they aren't.   Thus it shows that the homeomorphism
problem is in $\coRE$.  This is good enough by \Prop{p:recore} and
\Cor{c:plre}.

The algorithm also does not by itself determine whether the isometry
is orientation-preserving.   However, this is very little extra work.
Given $\eps \ll \delta$ and given $\eps$-nets $S_1 \subseteq N_1$ and $S_2
\subseteq N_2$, we can let $p_1,p_2,p_3,p_4$ be 4 points in $S_1$ that
lie in a ball of radius $\delta/2$ and that make an approximately regular
tetrahedron $\Delta$.   If $f:S_1 \to S_2$ is an approximate isometry, then
we can check whether $f$ flips over $\Delta$.  If no orientation-preserving
isometry exists, then when $\eps$ is small enough, either $f$ will cease
to exist or $\Delta$ will be flipped over.

We can use similar methods to find the mapping class group of a single $N$,
since by Mostow rigidity it is also the isometry group of $N$.  We assume
that $N$ has boundary, which is technically short of the full generality
of \Thm{th:hyphom}, but enough to prove \Thm{th:main}.  Just as with the
method to check whether and approximate $f$ preserves orientation, we can
when $\eps$ is small enough compute the effect of $f$ on $H_1(\del N)$,
which determines which isometry is close to $f$ (if any).
\end{proof}

\begin{proof}[Third proof] In this proof, we restrict attention to case 1 of
\Lem{l:hyprec} and thus work over the ring $\hQ$ of real algebraic numbers.
We assume real algebraic coordinates for $\H^3$ and for its isometry group
$\Isom^+(\H^3)$; for example we can take $\H^3$ to be the set of positive,
unit timelike vectors in $3+1$-dimensional Minkowski geometry, and we can
take $\Isom^+(\H^3) = \ISO(3,1)$.  We again assume that $N_1$ and $N_2$
are made from $N^*_1$ and $N^*_2$ using a common truncation area $2^{-n}$.

We assume geometric triangulations $\Theta^*_1$ and $\Theta^*_2$ of $N^*_1$
and $N^*_2$ with real algebraic descriptions.   Using these triangulations,
we can find finite, open coverings of $N_1$ and $N_2$ by metric balls
$B(p,\eps)$, where each point $p$ has a real algebraic position and
the common radius is (a) also real algebraic, and (b) less than half of
the injectivity radius of $N_1$ and $N_2$.  Then we can give each ball
the same algebraic coordinates as $\H^3$, and we can also calculate the
relative position of every pair of balls as some element in $\Isom^+(\H^3)$.
In other words, we obtain atlases of charts for $N_1$ and $N_2$ using the
$\Isom^+(\H^3)$ pseudogroup.  In fact, everything is constructed in the
subgroup and sub-pseudogroup with real algebraic matrix entries.

If there is an isometry between $N_1$ and $N_2$, then their atlases
combine into a larger atlas.  There are only finitely many possible
patterns of intersection between the balls of $N_1$ and the balls of $N_2$.
For each such pattern, we obtain a finite system of algebraic equalities
and inequalities, which says first that the intersection pattern is what is
promised, and second that the gluing maps between the atlases are consistent.
\Thm{th:tarski} then says that it is recursive to determine whether this
system of equations has a solution.  Since we work in the group
$\Isom^+(\H^3)$, we are looking only for orientation-preserving isometries.
\end{proof}

\section{Homeomorphism is in $\ER$}

We will use the basic fact that a finite composition of $\ER$ functions is
in $\ER$.  In other words, if an algorithm has a bounded number of stages
that expand its data by an exponential amount or otherwise by an $\ER$
amount, then it is still in $\ER$.

\subsection{The outer proof}

In this section we will prove \Thm{th:ermain}.  The proof is a combination of
the proof of \Thm{th:main} together with several computational improvements.
We summarize these computational improvements in this section by stated some
supporting theorems which we will prove ourselves (or prove by citation)
with two main supporting tools.  The first tool is normal surface theory,
which we can use to find essential spheres and tori and Seifert fibrations.
Note that Jaco, Letscher, and Rubinstein \cite{JLR:essential} sketched a
similar approach.  The second tool is an $\ER$ version of \Thm{th:tarski}
\cite{Grigoriev:algebra}, which we use to bound the complexity of a
geometric triangulation of a hyperbolic manifold, and the complexity of
recognizing small Seifert-fibered spaces.

\begin{theorem} It is in $\ER$ to find and triangulate the prime summands
$\{W\}$ of a closed, oriented 3-manifold $M$, to find and triangulate the
JSJ components $\{N\}$ within each prime summand $W$, to find their JSJ
graph $\Gamma$, and to recognize which components $N$ are Seifert-fibered
and find their fibrations.
\label{th:erjsj} \end{theorem}

We will prove most of \Thm{th:erjsj} in \Sec{s:normal} using normal surface
theory.  Note that Jaco, Letscher, and Rubinstein \cite{JLR:essential}
sketched ideas that are similar to our proof.

\Thm{th:erjsj} also has one lingering case which is more difficult.
Recall that a Seifert-fibered space is \emph{small} if it is non-Haken
(and therefore closed).

\begin{theorem} Recognizing small Seifert-fibered spaces is in $\ER$.
\label{th:ersmall} \end{theorem}

We will prove \Thm{th:ersmall} in \Sec{ss:small} using a combination of
normal surface theory and algebraic methods.  Note that Li \cite{Li:small}
shows that recognizing small Seifert-fibered spaces with infinite $\pi_1$
is recursive, and his algorithm should be elementary recursive.  However,
we will use a different approach for this part of the theorem.  Li also
addresses the finite $\pi_1$ case in two different ways.  Without assuming
geometrization (which was still open at the time), he cites work of
Rubinstein and Rannard-Rubinstein on small Seifert-fibered spaces.  He also
outlines a simplified argument for the finite $\pi_1$ caes that depends
on geometrization; we give a detailed argument which is in a similar spirit.

\begin{theorem} If a compact, oriented 3-manifold $N$ has a closed or
cusped hyperbolic structure, then it is in $\ER$ to find a geometric
triangulation and specify its geometric data with algebraic numbers.
The isomorphism and automorphism problems are also both in $\ER$.
\label{th:erhyp} \end{theorem}

We will prove \Thm{th:erhyp} in \Sec{ss:erhyp} using both
Mostow rigidity and methods from algebraic geometry.

\begin{proof}[Proof of \Thm{th:ermain}] We consider each stage of
the proof of \Thm{th:main} in turn.  The proof begins with a geometric
recognition of a single closed, oriented 3-manifold $M$ in \Thm{th:geomrec}.
This is not elementary recursive as described, but we can replace it with
\Thm{th:erjsj} to find the direct sum and JSJ decomposition.   We can
then apply \Thm{th:erhyp}, which is an $\ER$ version of \Lem{l:hyprec},
to calculate the hyperbolic structure of each hyperbolic summand $N$.
This includes a description of the Euclidean structure of each torus
component of $\del N$.

Finally given two closed, oriented 3-manifolds $M_1$, $M_2$, we first
decompose them into summands.  For each bijection among the summands, we want
to calculate $W_1 \congq W_2$ for each pair of matching summands.  This is
a calculation with JSJ graphs which is done in \Sec{ss:graph} to complete
the proof of \Thm{th:main}, and this part is already elementary recursive.
\end{proof}

\subsection{Normal surfaces}
\label{s:normal}

Let $M$ be a compact 3-manifold with triangulation $\Theta$.  Recall that
a \emph{normal surface} $S \subseteq M$ intersects each tetrahedron $\Delta
\in \Theta$ in 7 types of elementary disks, namely 4 types of triangles and
3 types of quadrilaterals.  The surface $S = S_v$ is given by a vector $v
\in \Z_{\ge 0}^{7t}$ that lists the number of each type of elementary disk.
If $v$ is such a vector, then $S_v$ is embedded (and uniquely defined)
provided that it only uses at most one type of quadrilateral in each
tetrahedron.  After specifying which type of quadrilateral is allowed in
each tetrahedron, the normal surface equations then have a polytopal cone
\[ C \subseteq \Z_{\ge 0}^{5t} \subseteq \Z_{\ge 0}^{7t} \]
of solutions.  We define a \emph{fundamental surface} $S_v$ is
one whose vector $v \in C$ is not the sum of two other solutions in $C$.
If $S_v$ is non-orientable, then $S_{2v}$ is its orientable double cover
and we call it fundamental as well.

\begin{lemma}[Haken] The number of elementary disks in a fundamental
surface in $M$ is bounded above by an exponential in $t$.  (Thus it is
elementary recursive.)
\label{l:fund} \end{lemma}

We can represent a normal surface $S$ by listing all triangles and
quadrilaterals in order in each tetrahedron $\Delta \in \Theta$.    It is
then easy to separate $S$ into connected components and calculate the
topology of each component.  This is exponentially inefficient compared
to algorithms such as Agol-Hass-Thurston \cite{AHT:genus}, but it has no
effect on whether the resulting algorithm is in $\ER$.

We define a \emph{complete set} of essential 2-spheres in a 3-manifold
$M$ is a collection $C$ such that cutting $M$ along each 2-sphere in $M$
and capping off the resulting boundary components produces irreducible
3-manifolds.  Likewise a complete set of essential disks is a collection
$C$ of properly embedded disks which are not boundary parallel, such that
the compression of every disk in $C$ renders $M$ boundary-incompressible.

\begin{theorem}[Jaco-Tollefson] Let $M$ be a compact, oriented, triangulated
3-manifold. Then:
\begin{enumerate}
\item $M$ has a collection of disjoint, fundamental surfaces which form
a complete set of essential 2-spheres \cite[Thm.~5.2]{JT:complete}.
\item If $M$ has no essential 2-spheres, then it has a collection of
disjoint, fundamental surfaces which form a complete set of essential
disks \cite[Thm.~6.2]{JT:complete}.
\item If $M$ has no essential 2-spheres or disks, and it has an essential
torus, then it has one which is fundamental \cite[Cor.~6.8]{JT:complete}.
\item If $M$ has no essential 2-spheres or disks, and it has an essential
annulus, then it has one which is fundamental \cite[Cor.~6.8]{JT:complete}.
\end{enumerate}
\label{th:jt} \end{theorem}

Actually Jaco and Tollefson show that each type of surface described in
\Thm{th:jt} is a vertex surface, which is a special case of a fundamental
surface.    Also case 1 of \Thm{th:jt} is stated for closed manifolds,
but the proof is the same for manifolds with boundary.  Finally,
given \Lem{l:fund}, the surfaces produces by \Thm{th:jt} all have
an elementary recursive bound on their size.

We will use the following variation of \Thm{th:jt}.

\begin{theorem}[Hass-K. \cite{K:sphere}] If $M$ is a closed, oriented
3-manifold with a triangulation $\Theta$, then it has a collection of
disjoint normal surfaces which form a complete set of essential spheres
and tori, such that the total number of elementary disks is bounded above
by an exponential in $t$.
\label{th:hk} \end{theorem}

Briefly, \Thm{th:hk} uses a generalization of the normal surface
equations which we call the \emph{disjoint normal surface equations}.
(They are similar to the crushed triangulation technique defined by Casson
\cite{JLR:essential}.)  They are equations for a normal surface $S$ which is
disjoint from a fixed normal surface $R \subseteq M$.  To prove \Thm{th:hk},
we find each surface $S$ or $T$ sequentially as a fundamental surface,
relative to the union of previous surfaces.

\begin{theorem}[Rubinstein \cite{Rubinstein:recog}, Thompson
\cite{Thompson:recog}] Recognizing the 3-sphere $S^3$ is in $\ER$.
\label{th:ersphere} \end{theorem}

The proof of \Thm{th:ersphere} uses a variant known as \emph{almost normal
surfaces} that are allowed one exceptional intersection with a tetrahedron
that is either an octagon, or a triangle and a quadrilateral with a
connecting annulus.  The original papers only claim a recursive algorithm,
but the algorithm is based on normal surface theory.  In fact, the proof also
uses disjoint normal surface equations.  Schleimer \cite{Schleimer:sphere}
also refines the Rubinstein-Thompson algorithm to show that 3-sphere
recognition is in the complexity class $\NP$, which is a much better bound
than just $\ER$.

\begin{proof}[Proof of \Thm{th:erjsj}] As a first step, we check whether $M
\cong S^3$ using \Thm{th:ersphere}.  If not, we search over collections $C$
of normal surfaces in $M$ with a suitable elementary recursive complexity
bound in order to find a set of surfaces that meets the conclusion of
case 1 of \Thm{th:jt}.  To test whether a given collection $C$ is one
that we want, we first calculate whether each surface in it is a 2-sphere.
Then we can cut along all of the spheres (and retriangulate) and cap them
to make a multiset of summands of $M$.  For each non-separating sphere,
we create a separate $S^2 \times S^1$ summand.  What remains is a putative
prime factorization $\{W\}$, but we must check whether the summands are
irreducible and not $S^3$.  We can use \Thm{th:ersphere} to check that no
summand $W$ is $S^3$.  If not, then we can again use case 1 of \Thm{th:jt}
to look for an essential 2-sphere in $W$, and again use \Thm{th:ersphere}
to check whether it is essential.

For each summand $W$, we similarly search for a collection $C$ that meets
the conclusion of \Thm{th:hk}.  We can check that each surface in $C$ is
a torus.  We can then cut $W$ along $C$ to make a putative decomposition of
$W$ into atoroidal components $\{Q\}$.  We have the following possibilities
for each component $Q$:
\begin{enumerate}
\item $Q$ has an essential disk, which necessarily cuts it into a ball.
In this case, $Q \cong S^1 \times D^2$ is a solid torus.
\item $Q$ has no essential disk, but it has a separating essential annulus
that cuts it into two solid tori.  In this case $Q$ fibers
over a disk with two exceptional fibers.
\item $Q$ does not have a separating essential annulus, but it does
have a non-separating essential annulus that cuts it into a solid torus.
In this case $Q$ fibers over an annulus or a M\"obius strip with at most
one exceptional fiber.
\item $Q$ has a separating annulus that cuts it into two thickened tori.
In this case $Q \cong S^1 \times F$, where $F$ is a pair of pants.
\item $Q$ has an essential torus, specifically an incompressible torus
which is not boundary-parallel.
\item $Q = W$ is closed and has no essential torus.  In this case $Q$
is either hyperbolic or small Seifert-fibered.
\item $Q$ has boundary, and it has no essential disk, annulus, or torus.
In this case, $Q$ is hyperbolic.
\end{enumerate}
To see that this is an exhaustive list, we recall that if $Q$ is
Seifert-fibered with boundary and is atoroidal, then with one exception
its orbifold base is planar, and the total number of boundary circles
plus exceptional fibers is at most three.  The only exception is that $Q
\cong K^2 \ltimes I$ has a M\"obius strip base; but it also has its other
fibration with Seifert data $\{D^2,(2,1),(2,1)\}$.

We claim that we can recognize each possibility for $Q$ by a bounded number
of applications of \Thm{th:jt}; the recognition algorithm is therefore
in $\ER$.  Indeed, each case reduces to earlier cases when a putative
essential surface is found.  The most subtle case is case 5, where we can
check whether a candidate torus in $Q$ is essential by checking that it is
not compressible (case 1) and does not bound a thickened torus (case 3).
Note that if $Q$ is a solid torus (case 1) or a thickened torus that is
not glued to itself (case 3), or if $Q$ has an essential torus (case 5),
then the collection $C$ in $W$ should be rejected.

In case 6, we use \Thm{th:ersmall} to determine if $W$ is small
Seifert-fibered, and if so, its homeomorphism type.

If $Q$ is a thickened torus which is glued to itself, then $W$ is a
torus bundle over a circle, and we can find its monodromy matrix $A \in
\SL(2,\Z)$ with a homology calculation.  We can solve the conjugacy problem
in $\SL(2,\Z)$ using the usual trick that $\PSL(2,\Z) \cong C_3 * C_2$,
and thus determine the homeomorphism type of $W$. (Note that $W$ may be
either Sol, Nil, or Euclidean.)

Otherwise we can piece together the JSJ decomposition $\{N\}$ of $W$
from the (often non-unique) atoroidal decomposition $\{Q\}$.  In each
remaining case $Q$ has either 0 Seifert fibrations (if it is hyperbolic),
or 2 fibrations (if it is $K^2 \ltimes I$), or 1 fibration (in all other
Seifert-fibered cases).  Using the recognition of $Q$, we can calculate
its Seifert data and express the fibration of each boundary torus $T
\subseteq \del Q$ by the corresponding line $L \subseteq H_1(T)$.  We can
then piece together adjacent fibered components to make JSJ components,
when the fibrations match.   The Seifert data produced in this manner is
not necessarily canonical, but canonicalizing it is straightforward.
\end{proof}

\subsection{Algebraic algorithms}
\label{ss:algebraic}

We list several complexity bounds concerning algebraic numbers and solutions
to algebraic equations.

\begin{theorem}[Collins, Monk, Vorobiev-Grigoriev, W\"uthrich
{\cite[Thm. 4]{Grigoriev:algebra}}] Suppose that a set $S \subseteq \R^n$ is
defined by a finite set of polynomial equalities and inequalities over $\Q$.
Then it is in $\ER$ to calculate a representative finite set $F \subseteq
\hQ^n$, with one point $p \in F$ in each connected component of $S$.
\label{th:ertarski} \end{theorem}

\Thm{th:ertarski} is of course an $\ER$ version of \Thm{th:tarski}.  In the
statement of the theorem, each element $\alpha \in \hQ$ described by its
minimal polynomial $a(x) \in \Z[x]$ and an isolating interval $\alpha \in
[b,c]$ that contains exactly one root of $a(x)$.

\begin{lemma} If $\alpha \in \hQ$ is a non-zero complex root of a polynomial
$a(x) \in \Z[x]$, then there is an $\ER$ upper bound on $|\alpha|$ and
$|\alpha^{-1}|$.
\label{l:bound} \end{lemma}

\begin{proof} Let $n = \deg a$ and write
\[ a(x) = a_n x^n + a_{n-1} x^{n-1} + \dots + a_1 x + a_0. \]
We can assume without loss of generality that $a_0 \ne 0$.  If $|\alpha|
> \sum_k a_k$, then us that $|a_n \alpha^n|$ is larger than the total norm
of all of the other terms, so by the triangle inequality, $a(\alpha) \ne 0$.
This establishes $\sum_k a_k$ as an upper bound on $|\alpha|$.  For the lower
bound, we can observe that $\beta = \alpha^{-1}$ is a root of the polynomial
\[ b(x) = a_0 x^n + a_1 x^{n-1} + \dots + a_{n-1} x + a_n = x^n a(x^{-1}). \]
We can thus repeat the argument.
\end{proof}

If $x_1,\ldots,x_n$ are algebraic numbers, then we say that an algebraic
number $z$ is an \emph{integral primitive element} if each $x_k = f_k(z)$
for some integer polynomial $f_k \in \Z[x]$.  It is a result of Galois
that every finite set of algebraic numbers has a primitive element; we
are interested in a computationally bounded version.

\begin{theorem}[Koiran {\cite[Thm.~4]{Koiran:hilbert}}] 
If $x_1,\ldots,x_n$ are algebraic numbers, then they have an integer
primitive element $z$ which can be computed in $\ER$, and such that
polynomials $f_k$ with $f_k(z) = x_k$ can also be computed in $\ER$.
\label{th:koiran} \end{theorem}

Koiran states the result in the form $x_k = f_k(z)/a_k$, where the
denominator $a_k$ is an elementary recursive integer.  However, it is not
difficult to modify $z$ to eliminate these denominators.  (Proof: Let $f
\in \Z[x]$ be the minimal polynomial of $z$ and let
\[ z_1 = \frac{z}{f(0)\prod_k a_k}. \]
Then both $1/a_k$ and $f_k(z)$ are expressible as integer polynomials
in $z_1$.  Therefore, so is their product.)

\begin{lemma} If $h \in \Z[x]$ is an integer polynomial, then there is
a prime $p$ that can be computed in $\ER$ such that $h(x)$ has a root
in $\Z/p$.
\label{l:ermod} \end{lemma}

\begin{proof} If $d = \deg h$, then $h$ attains the values $\pm 1$
at most $2d$ times.  Therefore there is an integer $a$ with $|a| \le d$
such that $h(a)$ is not $\pm 1$, and we can let $p$ be a prime divisor
of $h(a)$.
\end{proof}

Using the results so far in this section, we obtain an elementary recursive
version of Mal'cev's theorem, which says that finitely generated residually
linear groups are residually finite.  Our computational version 
requires a finitely present

\begin{theorem} Let $\Gamma$ be a finitely presented group, let $g \in
\Gamma \setminus \{1\}$ be a non-trivial element given by a word $w$
in the generators of $\Gamma$, and suppose that there is a representation
\[ \rho:\Gamma \to \GL(n,\C) \]
that distinguishes $g$ from the identity. Then $\Gamma$ admits a finite
representation
\[ \rho_p:\Gamma \to \GL(n,\Z/p) \]
that distinguishes $g$ from the identity, where $p$ is a prime number.
Moreover, we can find such a $p$ and $\rho_p$ in $\ER$ given the presentation
of $\Gamma$, the word $w$, and the integer $n$.
\label{th:ermalcev} \end{theorem}

\begin{proof} A function from the generators of $\Gamma$ to $n \times n$
matrices forms a representation
\[ \rho:\Gamma \to \GL(n,\C) \]
if and only if the matrix entries satisfy equations that come from the
relators of the presentation of $\Gamma$.  We also want $\rho(g) \ne I$.
To this end, we assume another matrix of variables $Y$ and impose the
condition
\[\mathrm{tr}(Y(\rho(g)-I)) = 1. \]
By hypothesis $\rho$ and $Y$ exist, and \Thm{th:ertarski} then produces
an algebraic, elementary recursive solution.  By \Thm{th:koiran}, the
matrix entries are generated by an integer primitive element $z$, and by
\Lem{l:ermod}, we can replace $z$ by a residue $\alpha \in \Z/p$ for some
prime $p$ computable in $\ER$.  We thus get a modular representation
\[ \rho_p:\Gamma \to \GL(n,\Z/p) \]
and a matrix $Y_p$ over $\Z/p$ again with the same properties.  Since $Y_p$
exists, $\rho_p(g)$ cannot be the identity.
\end{proof}

\begin{remark} Assuming the Generalized Riemann Hypothesis, Koiran's work
implies \Thm{th:ermalcev} with a much better bound, namely that $\log(p)$
can be bounded by a polynomial in $n$, the length of the presentation of
$\Gamma$, and the length of the word $w$.
\end{remark}

\subsection{The small Seifert-fibered case}
\label{ss:small}

In this section we will prove \Thm{th:ersmall}.

Let $N$ be a closed, oriented 3-manifold which has been recognized
as irreducible and atoroidal by part of the algorithm in the proof
of \Thm{th:erjsj}.   We want to distinguish between the case that $N$
is small Seifert-fibered and the case that $N$ is hyperbolic; and in the
former case, find its homeomorphism type.  By \Thm{th:ersphere}, we also
may as well assume that $N \ne\cong S^3$.  We divide the proof into two
cases, according to whether $\pi_1(N)$ is finite or infinite.

\begin{proposition} It is in $\ER$ to determine if $\pi_1(N)$ is finite
and compute its oriented homeomorphism type.
\end{proposition}

\begin{proof} We first calculate whether $N$ is a lens space.  We calculate
$H_1(N)$ by applying the Smith normal form algorithm to its chain complex.
If $H_1(N)$ is infinite, then $N$ is not small Seifert-fibered.  Otherwise
the cardinality of $H_1(N)$ is elementary recursive.  We can calculate
whether $N$ is a lens space by checking whether $H_1(N) \cong \Z/m$ is cyclic
and calculating whether its abelian universal cover $\tN$ is isomorphic to
$S^3$.  To determine the parameter $n$ in the homeomorphism type of $N \cong
L(m,n)$, we can calculate the Reidemeister torsion of the twisted homology
of $N$ over the ring $\Z[\zeta]$, where $\zeta$ is an $m$th root of unity.
This is a determinant calculation which is a priori elementary recursive.
Note that this torsion determines the oriented homeomorphism type of $N$.

If $N$ is spherical but not a lens space, then it has a finite covering space
of order at most 12 which is a lens space.  We thus obtain an elementary
recursive bound on the cardinality of $\pi_1(N)$.  Using a presentation of
$\pi_1(N)$ obtained from the triangulation of $N$, we can search exhaustively
among surjective homomorphisms $\phi:\pi_1(N) \to \Gamma$, where $\Gamma$ is
a finite candidate for $\pi_1(N)$.  For each such surjective homomorphism,
we can build the corresponding covering space $\tN$ and calculate whether
$\tN \cong S^3$.  If this happens, then we know that $N \cong S^3/\Gamma$
as unoriented 3-manifolds.

Finally in the spherical case, we want to pass from the unoriented to the
oriented homeomorphism type of $N$ when $N$ is spherical but not lens.
(Note that every such $N$ is chiral.)  As a first warm-up, recall that
we can distinguish the lens space $L(4,1)$ from its reverse $L(4,-1) =
L(4,3)$ by computing its Reidemeister torsion.  As a second warm-up, we
consider the simplest prism space $R(1,2)$ whose fundamental group $\Gamma =
\pi_1(R(1,2))$ is the quaternionic 8-element group.   We can build $R(1,2)$
as a coset space inside $\SU(2)$:
\[ \Gamma \subseteq \SU(2) \qquad R(1,2) \cong \SU(2)/\Gamma. \]
The group $\Gamma$ has four cyclic subgroups of order 4 which are
not conjugate in $\Gamma$ itself, but which are conjugate in $\SU(2)$.
Matching this calculation to the Seifert data, $R(1,2)$ has three double
covers which are all-oriented homeomorphic to $L(4,1)$.  We can thus
calculate the orientation of $N$ by calculating whether any double cover
is $L(4,1)$ or $L(4,3)$.

In general if $N$ is spherical but not lens, then the quotient group
$\pi_1(N)/Z(\pi_1(N))$ is either a dihedral group, which is the prism
case; or the isometry group of a Platonic solid: tetrahedral, octahedral,
or icosahedral.  In the case that $N \cong R(m,n)$ and $m$ is odd, as well
as in the Platonic cases, $\pi_1(N)$ has a unique subgroup isomorphic to
the quaternionic group.  Thus we can form the corresponding covering space
$\tN \cong R(1,2)$ and calculate its orientation as in the second warmup.
Mean while if $N \cong R(m,n)$ and $m$ is even, then the center $Z(\pi_1(N))
\cong \Z/(2m)$ has a unique subgroup of order 4.  Thus $\pi_1(N)$ has a
canonically chosen cyclic subgroup of order 4, and we can again form $\tN$
and calculate whether it is $L(4,1)$ or $L(4,3)$.
\end{proof}

To prove \Prop{p:sfsinf}, we will use a different combinatorial model
of Seifert-fibered spaces than the one in \Sec{sss:seifcomp}.  If $N$
is Seifert-fibered with base $F$, then we can consider a triangulation
$\Theta$ of $F$ with the orbifold points placed at the vertices.  For each
triangle $\Delta \in \Theta$, we make a solid torus $\Delta \times S^1$
which we interpret as a chart for a circle bundle with structure group
$S^1$.  Then we can construct $N$ with an atlas of charts of this type.
(It is an atlas with closed charts rather than open charts, but this
is valid in context.)  When two triangles $\Delta_1$ and $\Delta_2$
intersect in an edge, we glue the charts together with a transition map
\[ f_{12}:\Delta_1 \cap \Delta_2 \to S^1 \cong \R/\Z. \] We can assume
that each transition map $f_{12}$ is affine-linear if lifted to $\R$, so
that the value of $f_{12}$ lies in $\Q/\Z$, and its slop is also in $\Q$.
Morevoer, it is not hard to convert the Seifert data for $N$ into these
transition functions, for any triangulation of $F$.  Finally, note that
if $p \in F$ is an orbifold point of order $a$ and $\Delta \in \Theta$
is any triangle that has $p$ as a vertex, then the gluing maps between
charts glue the circle over $p$ in such a way that it shortens by a
factor of $a$ and becomes a singular fiber of $N$.

\begin{proposition} It is in $\ER$ to determine if $N$ is small
Seifert-fibered with infinite $\pi_1(N)$, and if so, compute its oriented
homeomorphism type.
\label{p:sfsinf} \end{proposition}

\begin{proof} If $\pi_1(N)$ is infinite and $N$ is small Seifert-fibered,
then the base $F$ of $N$ is a 2-sphere with orbifold points of order $a_1
\ge a_2 \ge a_3$ with
\[ \frac1{a_1} + \frac1{a_2} + \frac1{a_3} \le 1. \]
In the equality case $F$ is Euclidean and $N$ is either Euclidean or Nil;
otherwise $F$ is hyperbolic and the geometry of $N$ is either $\H^2
\times \R$ or $\widetilde{\Isom(\H^2)}$.  The orbifold fundamental
group $\pi_1(F)$ is a von Dyck group $D(a_1,a_2,a_3)$ which is the
orientation-preserving subgroup of index two in the corresponding triangle
group $\Delta(a_1,a_2,a_3)$ in either $\Isom(\E^2)$ or $\Isom(\H^2)$,
and $\pi_1(N)$ is a central extension of $D(a_1,a_2,a_3)$ by $\Z$.

\begin{figure} \begin{center}
\begin{tikzpicture}[semithick]
\draw (-1.71,-1.71) -- (1.71,1.71);
\draw (1.47,-1.47) -- (-1.47,1.47);
\draw[bend right=4] (-1.47,1.47) to (1.71,1.71);
\draw[bend right=4] (1.47,-1.47) to (-1.71,-1.71);
\draw[bend right=4] (1.71,1.71) to (1.47,-1.47);
\draw[bend right=4] (-1.71,-1.71) to (-1.47,1.47);
\draw[bend left=3] (1.47,-1.47) to (2.92,-.19);
\draw[bend left=3] (-1.47,1.47) to (-2.92,.19);
\draw[bend left=4] (2.92,-.19) to (1.71,1.71);
\draw[bend left=4] (-2.92,.19) to (-1.71,-1.71);
\draw[bend left=3] (2.92,-.19) to (4.31,.74);
\draw[bend left=5] (4.31,.74) to (1.71,1.71);
\draw[bend right=3] (1.71,1.71) to (3.72,2.14);
\draw[bend right=6] (3.72,2.14) to (4.31,.74);
\fill (0,0) circle (.07);
\fill (-1.71,-1.71) circle (.07);
\fill (1.71,1.71) circle (.07);
\fill (1.47,-1.47) circle (.07);
\fill (-1.47,1.47) circle (.07);
\fill (2.92,-.19) circle (.07);
\draw[bend right=3,blue,very thick] (-2.92,.19) to (4.31,.74);
\draw[bend right=6,blue,very thick] (-2.92,.19) to (3.72,2.14);
\draw[bend right=3,red,very thick] (-2.92,.19) to (-1.27,1.27); 
\draw[bend right=4,red,very thick] (-1.27,1.27) to (1.71,1.53); 
\draw[bend right=5,red,very thick] (1.71,1.53) to (4.31,.74); 
\draw[bend right=3,red,very thick] (-2.92,.19) to (-1.37,1.37); 
\draw[bend right=4,red,very thick] (-1.37,1.37) to (1.61,1.61); 
\draw[bend right=3,red,very thick] (1.61,1.61) to (3.72,2.14); 
\draw[red,fill=blue] (-2.92,.19) circle (.1);
\draw[red,fill=blue] (4.31,.74) circle (.1);
\draw[red,fill=blue] (3.72,2.14) circle (.1);
\draw (2.2,-1.3) node {$\Phi$};
\draw[blue] (2,.1) node {$\Theta$};
\end{tikzpicture}
\caption{Isotopy of a triangulation $\Theta$ of $\tF$ (initially in
blue) so that in its new position (in red), its edges are parallel to a
$\pi_1(N)$-invariant triangulation $\Phi$.}
\label{f:parallel} \end{center} \end{figure}
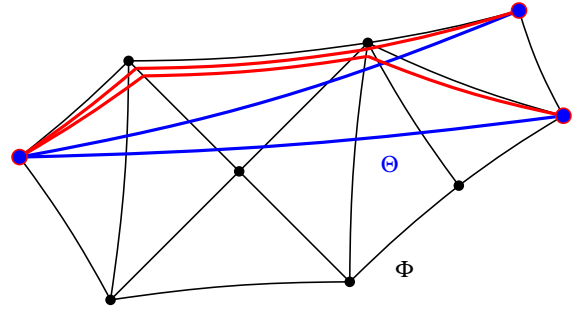

We first consider the case in which $F$ is Euclidean, which implies
that $(a_1,a_2,a_3)$ is either $(3,3,3)$, $(4,4,2)$, or $(6,3,2)$.
In this case there is a homomorphism 
\[ \phi:\pi_1(N) \to G, \]
where the target group is respectively the dihedral group $D_3$, $D_4$,
or $D_6$, such that the corresponding regular cover $\tN$ is a circle
bundle over a torus.  Given a putative choice for $G$ and $\phi$, we can
construct the regular cover $\tN$ and apply the large Seifert-fibered case
of \Thm{th:erjsj} to recognize it.  If $\tN$ is indeed a circle bundle
over a torus, then we know that $N$ must be small Seifert-fibered with a
Euclidean base, and the remaining question is to confirm that $G$ and $\phi$
were correctly chosen and thus compute the specific Seifert data of $N$.

Either $\tN$ is $S^1 \times S^1 \times S^1$ (so that $N$ itself
is Euclidean), or it is a circle bundle over $S^1 \times S^1$ with a
non-trivial Euler number and thus has Nil geometry.  We thus obtain
an explicit form of $\pi_1(\tN)$ which is either $\Z^3$ or a central
extension of $\Z^2$ by $\Z$.  At the same time, since the recognition of
the structure of $\tN$ is based on normal surface theory, it thus yields
a retriangulation of $\tN$ in $\ER$ from the one triangulation induced by
the input triangulation of $N$ to one that reveals the Seifert structure.
Therefore we obtain an explicit (and elementary recursive) description of
$\pi_1(N)$ as an extension of the finite group $G$ by $\pi_1(N)$.  We can
thus match $\pi_1(N)$ to the corresponding model Seifert-fibered space
to determine the unoriented homeomorphism type of $N$.  Finally we have
to calculate the oriented homeomorphism type.  In the Euclidean case, the
recognition of $\tN$ gives us an orientation of $\pi_1(\tN) \cong \Z^3$.
Similarly in the Nil case, when $\pi_1(\tN)$ is an extension of $\Z^2$
by $\Z$, the orientation of $\tN$ still gives an orientation of both the
center $\Z$ and the quotient $\Z^2$, up to switching both orientations.
This fixes an orientation of $N$ itself.

The argument when the base $F$ is hyperbolic is similar to the Euclidean
case but more complicated.  In this case, $\pi_1(N)$ has a non-trivial
homomorphism to $\Isom(\H^2)$, which in turn embeds in $\SL(2,\C)$,
so \Thm{th:ermalcev} tells us that $\pi_1(N)$ has a non-trivial finite
quotient $G$ which we can find in $\ER$ even if we do not know the linear
representation that explains that it exists.

We construct the finite cover $\tN$ of $N$ corresponding to the quotient
map $\phi:\pi_1(N) \to G$, and we apply part of \Thm{th:erjsj} to determine
if $\tN$ is large Seifert-fibered, and if so calculate its fibration and
its base $\tF$.  As before, we can first learn from this whether $N$ is
indeed small Seifert-fibered.  Second, as before this part of \Thm{th:erjsj}
gives us an $\ER$ retriangulation from the initial triangulation of $\tN$
as a covering space of $N$, to a triangulation that reflects its Seifert
fibration.   In particular, we obtain a triangulation $\Theta$ of $\tF$
together with an atlas of charts to describe $\tN$.  Third, we can choose
an orientation of $\tF$ and an orientation of the circle fibers so that
the two orientations together are consistent with the orientation of $\tN$
inherited from $N$.  Fourth, using the orbifold point orders of $\tF$ and
the cardinality of $G$, we obtain an $\ER$ upper bound on $a_1$, $a_2$,
and $a_3$.

For each candidate for $(a_1,a_2,a_3)$, the given representation of
$D(a_1,a_2,a_3)$ in $\Isom(\H^3)$ is rigid.  It is still rigid even as
a representation of $\pi_1(N)$, because any nearby representation must
still annihilate the kernel $Z(\pi_1(N))$.   Passing to the covering space
$\tF$, we obtain a preferred hyperbolic structure on $\tF$ and we can
realize $\Theta$ as a geometric triangulation.  (In two dimensions, every
triangulation of a hyperbolic surface is geometric.)  Now \Thm{th:ertarski},
combined with the fact that the retriangulation of $\tN$ is in $\ER$,
gives us an $\ER$ upper bound on the lengths of the edges of $\Theta$.
At the same time, we get a second (possibly generalized) triangulation
$\Phi$ of $\tF$ by tiling it by lifts of the triangular fundamental domain
of the triangle group $\Delta(a_1,a_2,a_3)$; this triangulation is both
geometric and $\pi_1(N)$-invariant.  As illustrated in \Fig{f:parallel},
we can isotop $\Theta$ so that its edges run parallel to $\Phi$.
(Since this isotopy cannot introduce crossings, and since two edges of
$\Theta$ might follow some of the same edges of $\Phi$ as illustrated,
they cannot land exactly onto these edges.)  This isotopy shows that
$\Phi$ and $\Theta$ have a mutual refinement that can be found in $\ER$.
Thus we can search over retriangulations of $\tF$ in $\ER$ until we find
one that is $\pi_1(N)$-invariant.  We can then use this to compute the
Seifert structure on $N$, moreover preserving the orientation information
inherited from the fibration of $\tN$.
\end{proof}

\subsection{The hyperbolic case}
\label{ss:erhyp}

To prove \Thm{th:erhyp}, we will need a quick mutual corollary of
\Thm{th:ertarski} and the proof of \Prop{p:refine}.

\begin{corollary} If $\Theta_1$ is a finite simplicial complex with
$n_1$ simplices (of arbitrary dimension) and $n_2 \ge n_1$, then it is
in $\ER$ to produce a complete list of geometric subdivisions $\Theta_2$
of $\Theta_1$ with $n_2$ simplices.
\label{c:errefine} \end{corollary}

\begin{proof}[Proof of \Thm{th:erhyp}] Let $\Theta$ be the input
triangulation of $N$ as a compact manifold, and let $\Theta^*$ be the
result of adding a cone to each component of $\del N$ to make a semi-ideal
combinatorial triangulation of $N^*$.  The manifold $N^*$ also has a
hyperbolic structure which we interpret as a separate manifold.  We rename
the hyperbolic version $X$ and assume a homeomorphism
\[f:N^* \to X. \]

\begin{figure}[htb]\begin{center}
\begin{tikzpicture}[scale=.8]
\useasboundingbox (-4.6,-5.15) rectangle (4.85,7.25);
\begin{scope}[shift={(-3,4.5)}]
\draw[fill=lightblue] (-1,0) -- (-1.5,1.5) -- (0,2) -- (1.5,1.5) -- (1,0);
\draw (-1,0) -- (1,0) -- (0,2) -- cycle;
\fill (-1,0) circle (2.5pt); \fill (1,0) circle (2.5pt);
\fill (-1.5,1.5) circle (2.5pt); \fill (1.5,1.5) circle (2.5pt);
\fill (0,2) circle (2.5pt);
\draw (0,-.75) node {$\Theta^*$};
\draw (0,2.1) node[anchor=south] {$1$};
\draw (-1.5,1.5) node[anchor=south east] {$2$};
\draw (-1,0) node[anchor=north east] {$3$};
\draw (1,0) node[anchor=north west] {$4$};
\draw (1.5,1.5) node[anchor=south west] {$5$};
\end{scope}
\draw[->] (-.75,5.5) -- (.75,5.5) node[midway,above] {$f$};
\begin{scope}[shift={(3,4.5)}]
\draw[fill=lightblue] (1,0) .. controls (1,3) and (-.5,3.5) .. (-1.5,1.5)
    -- (0,2) .. controls (-5,-2) and (4,-2) .. (1.5,1.5) -- (1,0);
\draw (-1,0) -- (1,0) -- (0,2) -- cycle;
\fill (-1,0) circle (2.5pt); \fill (1,0) circle (2.5pt);
\fill (-1.5,1.5) circle (2.5pt); \fill (1.5,1.5) circle (2.5pt);
\fill (0,2) circle (2.5pt);
\draw (0,2.1) node[anchor=south] {$1$};
\draw (-1.5,1.5) node[anchor=south east] {$5$};
\draw (-1,0) node[anchor=north east] {$3$};
\draw (1,0) node[anchor=north west] {$4$};
\draw (1.5,1.5) node[anchor=south west] {$2$};
\end{scope}
\begin{scope}[shift={(-3,0)}]
\draw[fill=lightblue] (-1,0) -- (-1.5,1.5) -- (0,2) -- (1.5,1.5) -- (1,0);
\draw (-1,0) -- (1,0) -- (0,2) -- cycle;
\fill (-1,0) circle (2.5pt); \fill (1,0) circle (2.5pt);
\fill (-1.5,1.5) circle (2.5pt); \fill (1.5,1.5) circle (2.5pt);
\fill (0,2) circle (2.5pt);
\draw (0,-.75) node {$\Theta^*$};
\draw (0,2.1) node[anchor=south] {$1$};
\draw (-1.5,1.5) node[anchor=south east] {$2$};
\draw (-1,0) node[anchor=north east] {$3$};
\draw (1,0) node[anchor=north west] {$4$};
\draw (1.5,1.5) node[anchor=south west] {$5$};
\end{scope}
\draw[->] (-.75,1) -- (.75,1) node[midway,above] {$g$};
\begin{scope}[shift={(3,0)}]
\fill[blue,opacity=.2] (-1,0) -- (1,0) -- (0,2) -- cycle;
\fill[blue,opacity=.2] (1,0) -- (-1.5,1.5) -- (0,2) -- cycle;
\fill[blue,opacity=.2] (0,2) -- (1.5,1.5) -- (-1,0) -- cycle;
\draw (-1,0) -- (1,0) -- (0,2) -- cycle;
\draw (1,0) -- (-1.5,1.5) -- (0,2) -- (1.5,1.5) -- (-1,0);
\fill (-1,0) circle (2.5pt); \fill (1,0) circle (2.5pt);
\fill (-1.5,1.5) circle (2.5pt); \fill (1.5,1.5) circle (2.5pt);
\fill (0,2) circle (2.5pt);
\draw (0,-.75) node {$\Lambda$};
\draw (0,2.1) node[anchor=south] {$1$};
\draw (-1.5,1.5) node[anchor=south east] {$5$};
\draw (-1,0) node[anchor=north east] {$3$};
\draw (1,0) node[anchor=north west] {$4$};
\draw (1.5,1.5) node[anchor=south west] {$2$};
\end{scope}
\begin{scope}[shift={(-3,-4.5)}]
\draw[fill=lightblue] (-1,0) -- (-1.5,1.5) -- (0,2) -- (1.5,1.5) -- (1,0);
\draw (-1,0) -- (1,0) -- (0,2) -- cycle;
\fill (-1,0) circle (2.5pt); \fill (1,0) circle (2.5pt);
\fill (-1.5,1.5) circle (2.5pt); \fill (1.5,1.5) circle (2.5pt);
\fill (0,2) circle (2.5pt);
\draw (0,-.75) node {$\Psi$};
\fill (.538,.923) circle (2.5pt); \fill (-.538,.923) circle (2.5pt);
\fill (0,.6) circle (2.5pt);
\fill (-1.2,.6) circle (2.5pt); \fill (-1.35,1.05) circle (2.5pt);
\fill (1.2,.6) circle (2.5pt); \fill (1.35,1.05) circle (2.5pt);
\draw (-1.35,1.05) -- (0,2) -- (-1.2,.6) -- (-.538,.923);
\draw (1.35,1.05) -- (0,2) -- (1.2,.6) -- (.538,.923);
\draw (0,.6) -- (0,2);
\draw (1,0) -- (-.538,.923); \draw (-1,0) -- (.538,.923);
\draw (0,2.1) node[anchor=south] {$1$};
\draw (-1.5,1.5) node[anchor=south east] {$2$};
\draw (-1,0) node[anchor=north east] {$3$};
\draw (1,0) node[anchor=north west] {$4$};
\draw (1.5,1.5) node[anchor=south west] {$5$};
\end{scope}
\draw[<-] (-.75,-3.5) -- (.75,-3.5) node[midway,above] {$g^{-1}$};
\begin{scope}[shift={(3,-4.5)}]
\fill[blue,opacity=.2] (-1,0) -- (1,0) -- (0,2) -- cycle;
\fill[blue,opacity=.2] (1,0) -- (-1.5,1.5) -- (0,2) -- cycle;
\fill[blue,opacity=.2] (0,2) -- (1.5,1.5) -- (-1,0) -- cycle;
\draw (-1,0) -- (1,0) -- (0,2) -- cycle;
\draw (1,0) -- (-1.5,1.5) -- (0,2) -- (1.5,1.5) -- (-1,0);
\draw (0,.6) -- (0,2);
\fill (-1,0) circle (2.5pt); \fill (1,0) circle (2.5pt);
\fill (-1.5,1.5) circle (2.5pt); \fill (1.5,1.5) circle (2.5pt);
\fill (0,2) circle (2.5pt);
\fill (.538,.923) circle (2.5pt); \fill (-.538,.923) circle (2.5pt);
\fill (0,.6) circle (2.5pt);
\draw (0,-.75) node {$\Phi$};
\draw (0,2.1) node[anchor=south] {$1$};
\draw (-1.5,1.5) node[anchor=south east] {$5$};
\draw (-1,0) node[anchor=north east] {$3$};
\draw (1,0) node[anchor=north west] {$4$};
\draw (1.5,1.5) node[anchor=south west] {$2$};
\end{scope}
\end{tikzpicture}\end{center}
\caption{We straighten $f$ to $g$, then subdivide the image, and
    finally subdivide the domain to make $g$ a simplicial map. In the
    proof the image is barycentrically subdivided, but any refinement
    which is a triangulation suffices.}
\label{f:erhyp} \end{figure}
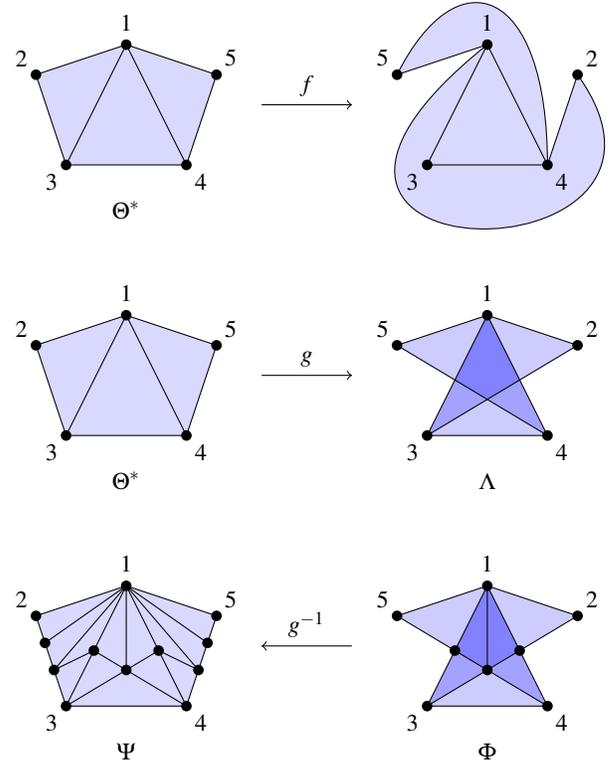

We fix the vertices of $N^*$ in the map $f$, and straighten all of
the tetrahedra, to make a map $g$ that represents $\Theta^*$ as a
self-intersecting geometric triangulation of $X$.   Since $g$ is homotopic
to $f$ (or properly homotopic of $N^*$ is not compact), it has (proper)
degree 1.  The simplices of $g$ may be flat or flipped over, but $g$
has degree 1 since it is homotopic to $f$.  The self-intersections of
$g(\Theta^*)$ yield a cellulation $\Lambda$ of $X$ with convex cells.
Thus $\Lambda$ has a barycentric subdivision $\Phi$ which is a geometric
triangulation of $X$.  Also let $\Psi = g^{-1}(\Phi)$.  Then $\Psi$ is a
refinement of the triangulation $\Theta^*$, and $g$ is now a simplicial
map from $\Psi$ to $\Phi$.  See \Fig{f:erhyp}.  (The figure uses a
simplicial refinement of the self-intersections which is simpler than
barycentric subdivision; this is not important for the proof.)

Now suppose that we do not know the hyperbolic structure of $N$, only that
it must have one because it prime, atoroidal, and acylindrical.  If we are
given $\Psi$ as a combinatorial refinement of the triangulation $\Theta^*$,
then we can search for $\Phi$ as a simplicial quotient of $\Psi$, such
that we can solve the hyperbolic gluing equations for $\Phi$ to recognize
it as a geometric triangulation of a hyperbolic manifold $X$.  We obtain
a candidate map $g:N^* \to X$.  If $g$ has degree 1, and there is also a
degree 1 map $h:X \to N^*$, then Mostow rigidity tells us that $g$ and $h$
are both homotopy equivalences and that $X$ and $N^*$ are homeomorphic.
(Note that there can be a degree 1 map in one direction between two
hyperbolic 3-manifolds that is not a homotopy equivalence \cite{BW:bundles},
even though this cannot happen in the case of hyperbolic surfaces.)  We can
search for $h$ by the same method of simplicial subdivision that we used
to find $g$.  This establishes an algorithm to calculate the hyperbolic
structure of $N^*$.

We claim that a modified version of this algorithm is in $\ER$.
We first consider $\ER$ candidates for the map $g$.  To do this, we
make a non-commutative cocycle $\alpha \in C^1(N;G)$ as in the proof of
\Cor{c:rigid}, where
\[ G = \Isom^+(\H^3) \cong \ISO(3,1),\]
and with the extra restriction that $\alpha$ is parabolic on each
component of $\partial N$.  These cocycle equations are algebraic, so
\Thm{th:ertarski} guarantees a representative set of solutions.  By Mostow
or Calabi-Weil rigidity, one of components of the solution space yields
a discrete homomorphism
\[ \rho:\pi_1(N) \to \Isom^+(\H^3) \]
that describes the hyperbolic geometry of $X$.  If we assign some point
$p \in \H^3$ to one of the vertices of $\Theta$, then in the closed case,
its orbit under $\alpha$ is in $\ER$ and can be extended on
each simplex of $\Theta$ to the map $g$.  In the cusped case, there are also
ideal vertices whose position on the sphere at infinity can be calculated
from $\alpha$ as well.

If $N$ has boundary, then we also want a truncated version of $\Theta^*$
which is larger than the original $\Theta$, and slightly different from
the horospheric truncation description in \Sec{s:hyphom}.  If that $\Delta
\in \Theta^*$ is semi-ideal, then let $p$ be its ideal vertex, let $F$
be the hyperplane containing the face of $\Delta$ opposite to $p$, and
let $F'$ be the hypersphere at distance $\log(2)$ from $F$ which is on the
same side as $p$.  Then we truncate $\Delta$ with $F'$ to make $\Delta'$;
or if $\Delta \in \Theta^*$ is a non-ideal tetrahedron, we let $\Delta'
= \Delta$.  We let $X' \subseteq X$ be the union of all $\Delta'$.
(In the closed case, we obtain $X' = X$.)  $X'$ can have a complicated
shape because the truncations are usually mismatched, but we can calculate
the positions of its vertices, and it is easy to confirm that it has at
least half of the volume of $X$.

Our algorithm does not know which cocycle $\alpha$ gives us a desired $g$
and we do not compute this directly.  Instead, we can calculate an $\ER$
bound for its data complexity, using the complexity bounds in the statement
of \Thm{th:ertarski}.

In particular, the existence of the map $g$ gives us $\ER$ bounds on the
parameters used in the third proof of \Lem{l:inject}.  Using \Lem{l:bound},
the existence of $g$ yields an $\ER$ upper bound on the diameter $\ell$
of $X'$ and then a lower bound on its injectivity radius $r$.

We can now follow the first proof of \Thm{th:hyphom}.  If $\Delta_1,
\Delta_2 \in \Theta^*$ are two tetrahedra, then the intersection complexity
of $g(\Delta_1)$ and $g(\Delta_2)$ is no worse than that of $g(\Delta'_1)$
and $g(\Delta'_2)$, and is bounded by an $\ER$ function of $\ell$ and $r$.
This yields an $\ER$ bound on the complexity of the refinements $\Phi$
and $\Psi$.  Recall that $\Psi$ is a refinement of $\Theta^*$, which is a
slightly modified version of the input description of $N$.  Having bounded
the complexity of $\Psi$, we can search for it using \Cor{c:errefine}
and solve for $\Phi$ and its geometry.  We can also discard $g$
if it does not have degree 1.

Thus far, the algorithm finds an $\ER$ collection of candidate maps $g:X^*
\to N$ of degree 1, where $N$ varies as well as $g$.  At least one of
these maps is a homotopy equivalence.   Instead of finding an inverse
$h$, We can repeat the algorithm to look for degree one maps among the
target manifolds $\{N\}$.   This induces a transitive relation among
these manifolds.  If $N$ is chosen at the top of this relation, then the
associated map $g:X^* \to N$ must be a homotopy equivalence.

To solve the isomorphism problem, we find geometric triangulations
$\Phi_1$ and $\Phi_2$ of the manifolds $N_1^*$ and $N_2^*$.  We can
again follow the first proof of \Thm{th:hyphom}, except now with an $\ER$
bound on the complexity of $\Phi_1$ and $\Phi_2$, and we can again
use \Cor{c:errefine}.  The same argument applies for the calculation of
the isometry group of a single $N^*$.
\end{proof}

\section{Open problems}
\label{s:final}

\Thm{th:ermain}, together with the fact that $\ER$ is a fairly
generous complexity class, suggests the following conjectures.

\begin{conjecture} If $M$ is a closed, Riemannian 3-manifold, then Ricci
flow with surgery on $M$ can be accurately simulated in $\ER$.
\label{c:ricci} \end{conjecture}

In other words, we conjecture that Perelman's proof of geometrization
can be placed in $\ER$.

\begin{conjecture} Every closed, hyperbolic manifold $N$ has a finite-sheeted
Haken covering which is computable in $\ER$.
\label{c:ervh} \end{conjecture}

In other words, we conjecture that the statement of the virtual Haken
conjecture, now the theorem of Agol et al \cite{Agol:haken}, can be placed
in $\ER$.   Maybe the known proof can be as well.

\begin{conjecture} Any two triangulations of a closed 3-manifold $M$
have a mutual refinement computable in $\ER$.
\label{c:erref} \end{conjecture}

\Conj{c:erref} does not follow from our proof of \Thm{th:ermain}, because
the algorithm in \Thm{th:erhyp} only establishes a simplicial homotopy
equivalence and then relies on Mostow rigidity.   However, the rest of the
proof of \Thm{th:ermain} uses a bounded number of normal surface dissections,
which does establish an $\ER$ mutual refinement according to the arguments
of Mijatovi{\'c} \cite{Mijatovic:s3,Mijatovic:seifert}.  Also, \Conj{c:ervh}
and the Haken case of \Conj{c:erref} would together imply the hyperbolic case
of \Conj{c:erref}, which would then imply the full conjecture.  Mijatovi{\'c}
\cite{Mijatovic:haken} also established that any two triangulations of a
fiber-free Haken 3-manifold have a primitive recursive mutual refinement.

Cases 3 and 4 of \Prop{p:er} are expected to be false for typical bounds
on complexity that are better than $\ER$.   Thus, in discussing further
improvements to \Thm{th:ermain}, we should consider qualitative complexity
classes, such as the famous $\NP$, rather than just bounds on execution time.
For one thing, $\ER$ is the union of an alternating, nested sequence of
time and space complexity classes, as follows:
\[ \ccP \subseteq \PSPACE \subseteq \EXP \subseteq
    \mathsf{EXPSPACE} \subseteq \mathsf{EEXP} \subseteq \cdots\;\;. \]
Here $\ccP$ is the set of decision problems that can be solved in
deterministic polynomial time; $\PSPACE$ is solvability in polynomial
space with unrestricted (but deterministic) computation time; $\EXP$ is
deterministic time $\exp(\poly(n))$; etc.  The author does not know where a
careful version of our proof of \Thm{th:ermain} would land in this hierarchy.


\providecommand{\bysame}{\leavevmode\hbox to3em{\hrulefill}\thinspace}
\providecommand{\MR}{\relax\ifhmode\unskip\space\fi MR }
\providecommand{\MRhref}[2]{%
  \href{http://www.ams.org/mathscinet-getitem?mr=#1}{#2}
}
\providecommand{\href}[2]{#2}
\providecommand{\eprint}{\begingroup \urlstyle{tt}\Url}

\typeout{get arXiv to do 4 passes: Label(s) may have changed. Rerun}

\end{document}